\documentclass[aop,MSNbibl,citesort,dvips]{arximspdf}
\usepackage{graphicx}

\doi{10.1214/11-AOP651}
\volume{39}
\issue{5}
\pubyear{2011}
\firstpage{1938}
\lastpage{1982}

\makeatletter

\def\bptnote#1{}

\newtheorem{theorem}{Theorem}
\newtheorem{lemma}[theorem]{Lemma}
\newtheorem{proposition}[theorem]{Proposition}
\newtheorem{corollary}[theorem]{Corollary}
\newtheorem{example}{Example}

\newtheorem{question}{Question}

\newcommand{\iint}{\int\!\!\int}
\newcommand{\iiint}{\int\!\!\int\!\!\int}

\newcommand{\Z}{{\mathbb Z}}
\newcommand{\R}{{\mathbb R}}
\newcommand{\N}{{\mathbb N}}
\newcommand{\X}{{\Omega}}
\newcommand{\w}{{\omega}}
\newcommand{\PP}{{\mathbb P}}
\newcommand{\E}{{\mathbb E}}
\newcommand{\borel}{{\mathcal{B}}}
\newcommand{\eqd}{\stackrel{d}{=}}
\newcommand{\F}{{\mathcal F}}
\newcommand{\powerset}{{\bolds{\mathcal P}}}
\newcommand{\vp}{{\mathcal{V}}}

\newcommand{\XX}{{\mathbb M}}
\newcommand{\XXd}{{\mathbb M'}}
\newcommand{\ww}{{\mu}}
\newcommand{\A}{{\mathcal A}}
\newcommand{\plusd}{+ d}
\newcommand{\thetaa}{\bolds{\sigma}}
\newcommand{\h}{{h'}}
\newcommand{\Phii}{{\Gamma}}
\newcommand{\standthin}{\phi^{\mathrm{ind}}}
\newcommand{\relstandthin}{F^{\mathrm{coin}}}
\newcommand{\gammathin}{\phi^{\mathrm{fin}}}
\newcommand{\poiseed}{\phi^{\mathrm{P}}}
\newcommand{\thmtwo}{{\Upsilon}}
\newcommand{\poiseedd}{\phi^{\mathrm{p}}}
\newcommand{\borelf}{{\mathfrak F}}
\newcommand{\borelbd}{{\mathcal B}^{*}}
\newcommand{\outside}{(\ssetglobe(X))^c}
\newcommand{\outsidemu}{(\ssetglobe(\mu))^c}
\newcommand{\lx}{\lambda}
\newcommand{\ly}{\lambda'}
\newcommand{\lyy}{\lambda- \lambda'}
\newcommand{\order}{\kappa}
\newcommand{\onee}{{\mathbf{1}}}
\newcommand{\bfz}{{\mathbf{0}}}
\newcommand{\ronn}{{|}}
\newcommand{\ronnn}{{|}}
\newcommand{\law}{{P}}

\newcommand{\emptypm}{0}
\newcommand{\modplus}{\oplus}

\newcommand{\gothh}[1]{{\mathcal #1}}
\newcommand{\goth}[1]{{\mathfrak #1}}

\newcommand{\les}{{\mathcal L}}

\newcommand{\sglobe}{\operatorname{Globes}^{1,2}}
\newcommand{\ssetglobe}{\Psi^{1,2}}
\newcommand{\oneglobe}{\operatorname{Globes}^1}
\newcommand{\twoglobe}{\operatorname{Globes}^2}
\newcommand{\onesetglobe}{\Psi^{1}}
\newcommand{\twosetglobe}{\Psi^{2}}

\newcommand{\outsideone}{(\onesetglobe(X))^c}
\newcommand{\outsidetwo}{(\twosetglobe(X))^c}
\newcommand{\outsideonemu}{(\onesetglobe(\mu))^c}
\newcommand{\outsidetwomu}{(\twosetglobe(\mu))^c}

\newcommand{\cglobe}{\mathbf{c}}
\newcommand{\ass}{\mathbf{U}}
\newcommand{\phithmtwo}{\Phi'}

\newcommand{\uniforma}{\mathrm{U}[0,1]}
\newcommand{\uniformb}{\mathrm{U}}

\newcommand{\ind}{{\mathbf1}}

\makeatother

\begin{document}
\begin{frontmatter}

\title{Poisson splitting by factors}
\runtitle{Poisson splitting by factors}

\begin{aug}
\author[A]{\fnms{Alexander E.} \snm{Holroyd}\thanksref{t1}\ead[label=e1]{holroyd@microsoft.com}\ead[label=u1,url]{http://research.microsoft.com/\textasciitilde holroyd}},
\author[B]{\fnms{Russell} \snm{Lyons}\corref{}\thanksref{t2}\ead[label=e2]{rdlyons@indiana.edu}\ead[label=u2,url]{http://mypage.iu.edu/\textasciitilde rdlyons}} and
\author[C]{\fnms{Terry} \snm{Soo}\thanksref{t3}\ead[label=e3]{tsoo@uvic.ca}\ead[label=u3,url]{www.math.uvic.ca/\textasciitilde tsoo}}
\runauthor{A. E. Holroyd, R. Lyons and T. Soo}
\affiliation{Microsoft Research, Indiana
University
and University of Victoria}
\address[A]{A. E. Holroyd\\
Microsoft Research\\
1 Microsoft Way\\
Redmond, Washington 98052\\
USA\\
\printead{e1}\\
\printead{u1}} 
\address[B]{R. Lyons \\
Department of Mathematics\\
Indiana University\\
Bloomington, Indiana 47405-5701\\
USA\\
\printead{e2}\\
\printead{u2}}
\address[C]{T. Soo \\
Department of Mathematics and Statistics\\
University of Victoria\\
PO BOX 3060 STN CSC\\
Victoria, BC V8W 3R4\\
Canada\\
\printead{e3}\\
\printead{u3}}
\end{aug}

\thankstext{t1}{Supported in part by Microsoft and NSERC.}
\thankstext{t2}{Supported in part by Microsoft and NSF Grant DMS-07-05518.}
\thankstext{t3}{Supported in part by NSERC.}

\received{\smonth{8} \syear{2009}}
\revised{\smonth{2} \syear{2011}}

%
\begin{abstract}
Given a homogeneous Poisson process on $\R^d$ with intensity $\lx$, we
prove that it is possible to partition the points into two sets, as
a~deterministic function of the process, and in an isometry-equivariant
way, so that each set of points forms a homogeneous Poisson process,
with any given pair of intensities summing to $\lx$. In particular,
this answers a question of Ball [\textit{Electron. Commun. Probab.}
\textbf{10} (2005) 60--69], who proved that in $d=1$,
the Poisson points may be similarly partitioned (via a
translation-equivariant function) so that \textit{one} set forms a Poisson
process of lower intensity, and asked whether the same is possible for
all $d$. We do not know whether it is possible similarly to \textit{add}
points (again chosen as a deterministic function of a Poisson process)
to obtain a Poisson process of higher intensity, but we prove that this
is not possible under an additional finitariness condition.
\end{abstract}

%

\vspace*{12pt}
\end{frontmatter}

\section{Introduction}\label{intro}

Let $\borel= \borel(\R^d)$ be the Borel $\sigma$-field on
$\R^d$. Let $\XX$ be the space of all Borel simple point
measures on $(\R^d, \borel)$, and let $\mathcal M$ be the
product $\sigma$-field on $\XX$ (we give detailed definitions
in Section \ref{examples}). Given an isometry~$\theta$ of
$\R^d$ and $\mu\in\XX$, we define $\theta(\ww)$ to be the
measure given by $\theta(\ww)(A) = \ww(\theta^{-1}(A))$ for all
$A \in\borel$. We say that a measurable mapping $\phi\dvtx\XX
\to\XX$ is \textit{isometry-equivariant} if
$\theta(\phi(\mu)) = \phi(\theta(\mu))$ for all $\mu\in\XX$
and for all isometries~$\theta$ of $\R^d$. Similarly\vspace*{1pt} we say
that $\phi$ is \textit{translation-equivariant} if it commutes
with all translations of $\R^d$. We define a partial order
$\leq$ on~$\XX$ via $\ww_1 \leq\ww_2 $ if and only if
$\ww_1(A) \leq\ww_2(A)$ for all $A \in\borel$. We say that
a~mapping~$\phi$ is \textit{monotone} if either $\phi(\mu) \leq
\mu$ for all $\mu\in\XX$, or $\mu\leq\phi(\mu)$ for all
$\mu\in\XX$.

Our main result is the following.
\begin{theorem}
\label{result} For all $ d \geq1$ and for all $\lx> \ly> 0$,
there exists a monotone isometry-equivariant mapping $\phi\dvtx\XX
\to\XX$ such that if $X$ is a homogeneous Poisson point
process on $\R^d$ with intensity $\lx$, then $\phi(X)$ and $X -
\phi(X)$ are homogeneous Poisson point processes on $\R^d$ with
intensities $\ly$ and $\lx- \ly$, respectively.
\end{theorem}

In other words, Theorem \ref{result} states that the points of
a Poisson process may be colored red and blue, in a
deterministic isometry-equivariant way, so that both the red
process and the blue process are Poisson processes.
Ball~\cite{ball} proved that in the case $d=1$, for all $\lx> \ly
> 0$, there exists a monotone translation-equivariant mapping
$\phi\dvtx\XX\to\XX$ such that if $X$ is a Poisson point process
with intensity $\lx$, then $\phi(X)$ is a homogeneous Poisson
point process with intensity $\ly$ (in other words, the Poisson
process may be ``thinned'' in a deterministic
translation-equivariant way). Ball asked whether the same is
possible in higher dimensions, and also whether the condition
of translation-equivariance can be strengthened to
isometry-equivariance. Theorem \ref{result} answers both
questions affirmatively, and also provides the additional
property that $X-\phi(X)$ is a Poisson process. Evans
\cite{Evans} recently proved that Poisson processes cannot be
thinned in an equivariant way with respect to any affine
measure-preserving group that is strictly larger than the
isometry group.

If all considerations of monotonicity are dropped, then the
following result of Ornstein and Weiss applies, even without
the restriction that $\lx>\ly$.
\begin{theorem}[(Ornstein and Weiss)]
\label{resulttwo} For all $ d \geq1$ and all $\lx, \ly
\in(0,\infty)$, there exists an isometry-equivariant mapping
$\phi\dvtx\XX\to\XX$ such that if $X$ is a~homogeneous Poisson
point process on $\R^d$ with intensity $\lx$, then $\phi(X)$ is
a~homogeneous Poisson point process on $\R^d$ with intensity
$\ly$.
\end{theorem}

Ornstein and Weiss \cite{MR910005} proved Theorem
\ref{resulttwo} as part of a much more general theory. In
particular, they proved the existence of an isomorphism,
whereas Theorem~\ref{resulttwo} asserts the existence only
of a homomorphism. The tools we develop to prove Theorem
\ref{result} allow us to give an alternative proof of
Theorem~\ref{resulttwo}. The map we construct is explicit, and it
satisfies an additional continuity property (see
Theorem~\ref{finitary} below). In addition, the map we construct is
\textit{source-universal}; that is, in Theorem \ref{resulttwo} the map $\phi
$ does not have to depend on the intensity of $X$. When $\lambda
'>\lambda$, we do not know
whether the condition of monotonicity can be added to Theorem
\ref{resulttwo} (in other words, whether a Poisson process can
be deterministically ``thickened'').
\begin{question}
\label{thick} Let $ d \geq1$ and let $\ly> \lx> 0$. Does
there exists a monotone isometry-equivariant $\phi\dvtx\XX\to
\XX$ such that if $X$ is a homogeneous Poisson point process on
$\R^d$ with intensity $\lx$, then $\phi(X)$ is a homogeneous
Poisson point process on $\R^d$ with intensity $\ly?$
\end{question}

However, we can prove that the answer to Question \ref{thick}
becomes \textit{no} when $\phi$ is required to satisfy the
following additional condition. For $\ww\in\XX$, we define
the restriction of $\ww$ to a set $A \in\borel$ via: $\ww
\ronn_A(\cdot):= \ww(\cdot\cap A)$ (so $\ww\ronn_A \in\XX$).
Let \mbox{$\| \cdot\|$} be the Euclidean norm on $\R^d$. The open
\textit{ball} of radius $r$ centered at $x$ is denoted by
$B(x,r): = \{ y \dvtx\|x-y\| < r \}$. Let $X$ be a Poisson point
process on $\R^d$ with law $\law$. We say that a
translation-equivariant measurable mapping $\phi\dvtx\XX\to\XX$
is \textit{strongly finitary} with respect to $\law$ if, for
$\law$-a.e. $\mu\in\XX$, there exists a positive real number
$n=n(\mu)$ such that for $\law$-a.e. $\mu' \in\XX$, we have
$\phi(\mu) \ronnn_{B(\bfz,1)} = \phi(\mu')
\ronnn_{B(\bfz,1)}$ whenever $\mu\ronnn_{B(\bfz,n)} =
\mu'\ronnn_{B(\bfz,n)}$. [In other words, the restriction of
$\phi(\mu)$ to the unit ball is determined by the restriction
of $\mu$ to a~larger ball, of random but finite radius.] With
the addition of this condition, we can answer Question
\ref{thick} in the negative, even if we drop the condition of
isometry-equivariance.
\begin{theorem}
\label{nothick} Let $d \geq1$ and $\ly> \lx> 0$. Let $X$
be a homogeneous Poisson point process on $\R^d$ with intensity
$\lx$ and law $\law$. There does not exist a~translation-equivariant monotone measurable mapping $\phi\dvtx\XX
\to\XX$ such that $\phi(X)$ is a homogeneous Poisson point
process on $\R^d$ with intensity $\ly$, and $\phi$ is strongly
finitary with respect to $\law$.
\end{theorem}

In fact, our proof of Theorem \ref{nothick} will not use the
assumption of translation-equivariance either, so we
actually prove the stronger statement that no mapping $\phi$
satisfying the other conditions can have have the property that
the restriction of $\phi(\mu)$ to the unit ball is determined
by the restriction of $\mu$ to a larger random ball, as defined
above.

In Section \ref{finitarysection}, we shall show that
the mappings that we produce to prove Theorems \ref{result}
and \ref{resulttwo} are strongly finitary. The mapping
produced in \cite{ball} is also strongly finitary.
\begin{theorem}
\label{finitary} Theorems \ref{result} and \ref{resulttwo} hold
even with the further requirement that the isometry-equivariant
mapping $\phi$ be strongly finitary with respect to~$\law$,
where~$\law$ is the law of $X$.
\end{theorem}

Sometimes deterministic translation-equivariant maps like the
ones of Theorems \ref{result} and \ref{resulttwo} are
called \textit{factors}. Factors are of basic importance in
ergodic theory and continue to play a central role in
applications of ergodic theory to combinatorics. The
combinatorial and probabilistic aspects of factors themselves
have received attention in recent years as well. It turns out
that factors are intimately related to Palm theory and
shift-coupling. For more information, see
\cite{Thorissontrv,last,Extra-Heads} and \cite{MR1741181}.
\textit{Factor graphs} of point processes have also received
considerable attention (see \cite{MR2044812,MR1961286,MR2081459}).
Following \cite{MR1961286}, a~\textit{factor graph} of a point process $X$ is a graph whose
vertices are the points of~$X$ and whose edges are obtained as
a deterministic translation-equivariant function of~$X$. An
important special case of a factor graph is a
translation-equivariant matching (see \cite{random} for some
striking results on this topic). Finally, we refer interested readers
to \cite{MR910005} for very general results
regarding factors of Poisson processes and the well-studied
isomorphism problem.

One can ask questions similar to ours about factors in a
discrete setting. Translation-equivariant matchings of i.i.d.
coin flips on $\Z^d$ are considered in~\cite{soo-2006} and~\cite{timarb}. Much is known about factors of Bernoulli
shifts on $\Z$ (e.g., see the monograph of Ornstein~\cite{MR0447525}). In particular, it is a classical result of
Sinai~\cite{MR0161960} that if $B(p)$ and $B(q)$ are
Bernoulli shifts on $\{0,1,\ldots, d-1 \}^{\Z}$ (i.e., i.i.d.
$\{ 0,1,\ldots, d-1\}$-valued sequences with laws $p$ and $q$),
and the entropy of $p$ is strictly greater than the entropy of
$q$, then there is a factor from $B(p)$ to $B(q)$. Recently,
Ball \cite{balltwo} proved that if the entropy of $p$ is
strictly greater than the entropy of $q$, and $p$
stochastically dominates $q$, then in the special case $d=2$,
there is a factor map $\phi$ from $B(p)$ to $B(q)$ that is
monotone [i.e., $\phi(x)_i \leq x_i$ for almost all $x \in
\{0,1,\ldots, d-1\}^{\Z}$ and all $i \in\Z$].

The factor\vspace*{1pt} map $\phi$ given in \cite{balltwo} is also
\textit{finitary}; that is, $\phi$ is continuous on a set of measure
one, when $\{0,1,\ldots, d-1\}^{\Z}$ is endowed with the product
topology. Keane
and Smorodinsky improved on results of Ornstein by producing
explicit finitary factors between Bernoulli shifts. We refer
the interested reader to the original papers of Keane and
Smorodinsky \cite{keanea,keaneb} and the recent survey
article on finitary codes by Serafin \cite{MR2306207}.

Finally, we also mention the work of Angel, Holroyd
and Soo \cite{AHS} concerning monotone deterministic functions
of Poisson point processes on finite volumes. In particular,
if $\lambda> \lambda'$, and $X$ is a Poisson point process of
intensity $\lambda$ on $[0,1]$, that article provides a
necessary and sufficient condition on $(\lambda, \lambda')$ for
the existence of a monotone deterministic map $\phi\dvtx\XX\to
\XX$ such that $\phi(X)$ is a Poisson point process on $[0,1]$
of intensity $\lambda'$.

\section{Some remarks about the proofs}
\label{examples} We next motivate the proofs of Theorems
\ref{result} and \ref{resulttwo} via some simple examples
of mappings $\phi\dvtx\XX\to\XX$ having some of the required
properties. The proof of Theorem \ref{nothick} is much
shorter and is treated in Section \ref{sectionthick}. Of
course, one of the requirements of $\phi$ is that it be
measurable. All the maps we define will clearly be measurable;
we provide the formal definition of the $\sigma$-field for
$\XX$ below.

\subsection*{Measurability}
The $\sigma$-field ${\mathcal M}$ of subsets of $\XX$ is
defined in the following way. Let $\N= \{ 0,1,2, \ldots \}$ be the
natural numbers, $\Z^{+} = \{ 1, 2, 3, \ldots \}$ be the positive
integers and $\bar{\N}$ be $\N\cup\{ \infty \}$. For
$B \in\borel$, the projection map $p_B\dvtx\XX\to\bar{\N}$ is
defined by $p_B(\mu) = \mu(B)$, for all $\mu\in\XX$. We let
${\mathcal M}$ be the smallest $\sigma$-field such that all the
projection maps are measurable.\vspace*{12pt}

Note that throughout this paper, the only laws we consider on
$\mathcal M$ will be homogeneous Poisson point processes on
$\R^d$ and their restrictions to subsets of $\R^d$. We say
that $U$ is a $\uniforma$ random variable if it is uniformly
distributed in $[0,1]$. Let $\les$
denote Lebesgue measure. Similarly, we say that~$V$ is a
$\uniformb[B]$ random variable if it is uniformly distributed in some
Borel set~$B$ with finite nonzero Lebesgue measure; that is, $\PP(V \in
\cdot) = {\les(\cdot\cap B) }/{\les(B)} $.
In the next examples and throughout
this paper, we shall assert that certain random variables can
be expressed as functions of $\uniforma$ random variables. This
can be justified by appealing to the Borel isomorphism theorem
\cite{borel}, Theorem~3.4.24. However, very often we need only the
following two results, which are consequences of the Borel
isomorphism theorem. Because of the need for isometry-equivariance in
our constructions, we shall often need to be
rather explicit about such functions.
\begin{lemma}[(Reproduction)]
\label{reprod} There exist measurable deterministic functions
$\{ g_i \}_{i \in\N}$, where $g_i\dvtx[0,1] \to[0,1]$, such that
if U is a $\uniforma$ random variable, then $\{ g_i(U) \}_{
i \in\N} $ is a sequence of i.i.d. $\uniforma$ random
variables.
\end{lemma}

For an explicit proof, see, for example, \cite{MR1876169}, Lemma 3.21.

Let $X$ be a Poisson process of intensity $\lx$ on $\R^d$. We say
that $Z$ is a~Poisson process of intensity $\lx$ on a set $A$ if $Z
\eqd X \ronn_A$.
\begin{lemma}[(Coupling)]
\label{poiseed}
Let $\ly> 0$. There exists\vspace*{1pt} a collection of
measurable mappings $\poiseed= \{ \poiseed_A \}_{A \in
\borel}$, where for each $A \in\borel$, the map $\poiseed_A =
\poiseed_{(A, \ly)}\dvtx[0,1]\to\XX$ is such that if $U$ is a
$\uniforma$ random variable, then $\poiseed_A(U)$ is a Poisson
point process on $A$ with intensity $\ly$.
\end{lemma}
\begin{pf}
By the Borel isomorphism theorem there exists a measurable function
$g\dvtx[0,1] \to\XX$ such that if $U$ is a $\uniforma$ random variable,
then $g(U)$ is a Poisson point process on $\R^d$ with intensity $\ly
$. Set $\poiseed_A(U) := g(U) \ronnn_A$.
\end{pf}
\begin{example}[(A $\Z^d$-translation-equivariant mapping
between Poisson point processes of arbitrary intensities)]
\label{firstfirsteg}
Let $\ly> 0$. Let $X$ be a
Poisson point process on $\R^d$ with positive intensity and law
$\law$. Let $C_\bfz$ be a cube of side-length~1 containing the
origin $\bfz\in\Z^d$, and let $C_i := C_\bfz+ i$ for $i \in
\Z^d$. Assume that~$C_\bfz$ is such that the collection
${\mathcal P} = \{ C_i \}_{i \in\Z^d}$ is a partition of
$\R^d$. The mapping~$\phi$ will be defined by specifying
$\phi(\cdot) \ronn_C$ for all $C \in{\mathcal P}$. We shall
define $\phi$ only off a~$\law$-null set; it is not difficult
to extend $\phi$ to all of $\XX$ so that it still commutes with
all translations of $\Z^d$. Let $g\dvtx[0,1] \to\XX$ be a
measurable function such that if $U$ is a $\uniforma$ random
variable, then $g(U)$ is a Poisson process on $C_\bfz$ with
intensity $\ly$. We shall define a measurable map $h\dvtx\XX\to
[0,1]^{\Z^d}$ with the following properties: $h(X)$ is a
collection of i.i.d. $\uniforma$ random variables, and for all
translations $\theta$ of $\Z^d$ we have $h(\theta(X))_i =
h(X)_{\theta^{-1}(i)}$ for all $i \in\Z^d$. For all $i \in\Z^d$,
let $\theta_i(x) = x+i$ for all $x \in\R^d$. Given the
mapping~$h$, it easy to see that by taking
\[
\phi(X)
\ronn_{C_i}:= \theta_{i} ( g(h(X)_i) )
\]
for all $i \in\Z^d$,
we have that $\phi$ commutes with translations of $\Z^d$ and
that $\phi(X)$ is a Poisson point process on $\R^d$ with
intensity $\ly$. It remains to define $h$.

If $X(C) = 1$, then we say that $C$ is \textup{special}. Let
$K(i)$ be the index of the first special cube to the right of
cube $i$; that is, $K(i) = i + (n, 0,\ldots, 0)$ where $n=n(i)$
is the smallest nonnegative integer such that $C_{i + (n, 0,\ldots,
0)}$ is special. Note that $\law$-a.s. $K$ is well
defined. For each special cube $C_i$, let~$z(i)$ be the
unique point $x \in C_i$ such that $X(\{ x \}) = 1$. Since $X$
is a Poisson point process, the random variables $\{ X \ronnn_{C_i}
\}_{i \in\Z^d}$ are independent, and also
conditional on the event that $C_i$ is special, $z(i)$ is a
$\uniformb[C_i]$ random variable. Let $f\dvtx C_\bfz\to
[0,1]^{\N}$ be a measurable function such that if $V$ is a
$\uniformb[C_\bfz]$ random variable, then $f(V)$ is a sequence
of i.i.d. $\uniforma$ random variables. For all $i \in\Z^d$,
let
\[
h(X)_{i}:= f\bigl( z(K(i)) - K(i) \bigr)_{n(i)}.
\]
It is
easy to verify that $h$ satisfies the required properties.
\end{example}

Let us remark that in Example \ref{firstfirsteg}, the map $\phi$ does
not depend on the intensity of $X$ and thus is source-universal. The
most important fact we used was that if $X$ is a Poisson process, then
conditional on the fact that it
has one point in $A$, the location of that point is a
$\uniformb[A]$ random variable. This elementary fact is true for any
Poisson process of positive intensity and will
often be useful. We shall appeal to it again in the next
example and in the proofs of Theorems \ref{result} and
\ref{resulttwo}. We refer the reader to \cite{kingman} or Theorem~1.2.1 of \cite{MR1199815}
for background and state
a slightly more general result in the lemma below.\looseness=-1
\begin{lemma}
\label{char} Let $X$ be a Poisson point process on $\R^d$ with
intensity $\lx$. Let $A \in\borel$ be a Borel set with positive finite
Lebesgue measure. Let $K$ be a Poisson random variable with
mean $\lx\les(A)$. Let $\{ V_i \}_{i \in\N}$ be a sequence of
i.i.d. $\uniformb[A]$ random variables that are independent of
$K$. Then $X \ronn_A$ has the same law as $Z:= \sum_{i=1} ^{K}
\delta_{V_i}$.
\end{lemma}

A central requirement in Theorems \ref{result},
\ref{resulttwo} and \ref{nothick} is that $\phi$ be a
deterministic function of $X$. The mapping\vspace*{1pt} in Example
\ref{firstfirsteg} is a deterministic function of~$X$ and
commutes with all translations of $\Z^d$. Given a
$\uniformb[C_\bfz]$ random variable~$V$, \mbox{independent} of $X$, we
can modify Example \ref{firstfirsteg} by starting with a~randomly shifted partition $\{ C_{i} + V \}_{i \in\Z^d}$ of
$\R^d$ and obtain a mapping $\Phi$ that is a~function of $X$
and $V$. As a result of starting with a randomly shifted
partition, the joint distribution of $(X, \Phi(X,V))$ is fully
translation-invariant. However,~$\Phi$ is no longer a
deterministic function of $X$.

Instead of using the lattice $\Z^d$, we shall use randomness
from the process to define a partition of $\R^d$. It is
straightforward to do this in an isometry-equivariant way. The
difficulty lies in choosing a partition that avoids potential
dependency problems.

We now turn our attention to Theorem \ref{result}. Let $\lx
> \ly> 0$. It is nontrivial to show that there exists a (not
necessarily translation-equivariant) monotone mapping which
maps a Poisson point process of intensity $\lx$ to a Poisson
point process of intensity $\ly$. In Example
\ref{firstfirsteg}, we asserted the existence of a certain
coupling between uniform random variables and Poisson point
processes via a measurable function $g\dvtx[0,1] \to\XX$ such that
whenever $U$ is a $\uniforma$ random variable, $g(U)$~is a
Poisson point process. Due to the monotonicity requirement in
Theorem \ref{result}, we require a more specialized coupling.

An important tool in the proof of Theorem \ref{result} will
be Proposition \ref{specquant} below, which is motivated by
one of the key ideas from Lemma 3.1 of \cite{ball}.
Proposition~\ref{specquant} provides a coupling between a Poisson point
process $X$ in a finite volume and another, $Y$, of lower
intensity, such that $Y \leq X$ and the process $X-Y$ is also a
Poisson point process. The process $Y$ is not a deterministic
function of $X$, but the coupling has certain other useful
properties.

Throughout this paper, it will be convenient to encode
randomness as a~function of $\uniforma$ random variables, as
was done repeatedly in Example~\ref{firstfirsteg}. For any
point process $Z$, the \textit{support} of $Z$ is the random set
\[
[Z] := \bigl\{ x \in\R^d \dvtx Z(\{ {x} \}) = 1 \bigr\}.
\]
Elements of
$[Z]$ are called \textit{$Z$-points}. We call a mapping
$\Phi\dvtx\XX\times[0,1] \to\XX$ a \textit{splitting} if
$\Phi(\mu, u) \leq\mu$ for all $(\mu, u) \in\XX\times[0,1]$,
and if for some $\lx> \ly$ we have that $\Phi(X, U)$ and $X - \Phi
(X, U)$ are Poisson point
processes with intensities $\ly$ and $\lx- \ly$, respectively,
whenever $X$ is a Poisson point process of intensity $\lx$, and $U$ is a
$\uniforma$ random variable independent of $X$. For example,
consider the coupling between a Poisson point process $X$ on
$\R^d$ of intensity $\lx$ and another, $Y$, of lower intensity
$\ly$, that is given by coloring the points of $X$
independently of each other red or blue with
probabilities~$\frac{\ly}{\lx}$ and $1-\frac{\ly}{\lx}$ and then taking the
red points to be the set of $Y$-points. It is easy to see that
both $Y$ (the red points) and $X-Y$ (the blue points) are
independent Poisson point processes on $\R^d$ with intensities
$\ly$ and $\lx-\ly$. This elementary result is sometimes
referred to as the \textit{coloring theorem} \cite{kingman}
and this coupling can be expressed as a splitting since all the
required coin-flips can be encoded as a function of a single
$\uniforma$ random variable. We shall revisit this elementary
coupling in more detail in Section \ref{selectionrule}. The
coupling given by Proposition \ref{specquant} below is also a
splitting.
\begin{proposition}[(Splitting on finite volumes)]
\label{specquant}
Let $\lx> \ly> 0$. There exists a finite constant $K= K(\lx,
\ly)$ and a family $\gammathin$ of measurable\vspace*{1pt}
mappings~$\gammathin_A$ so that for each $A \in\borel$ with finite
Lebesgue measure larger than~$K$, the map $\gammathin_A =
\gammathin_{(A,\lx,\ly)}\dvtx\XX\times[0,1] \to\XX$ has the
following properties:
\begin{longlist}[(a)]
\item[(a)]\hypertarget{monotonedist}
The map $\gammathin_A$ is monotone;
that is, $\gammathin_A(\mu, u) \leq\mu$ for all $(\mu, u) \in
\XX\times[0,1]$.
\item[(b)]
\hypertarget{monoquant} For all $(\mu, u) \in\XX\times[0,1]$, we
have $\gammathin_A(\mu, u) = \gammathin_A(\mu\ronn_A, u)$.
\item[(c)]
\hypertarget{dist} If $X$ is a homogeneous Poisson point process on
$\R^d$ with intensity~$\lx$, and $U$ is a $\uniforma$ random
variable independent of $X$, then
$\gammathin_A(X \ronn_A, U)$ is a Poisson point process of
intensity $\ly$ on $A$, and $X\ronn_A - \gammathin_A(X
\ronn_A, U)$ is a~Poisson point process of intensity $\lx-\ly$ on $A$.
\item[(d)]\hypertarget{interest} For all $(\mu, u) \in\XX\times[0,1]$,
if $\mu(A) = 1$, then $\gammathin_A(\mu,u) = 0$, while if
$\mu(A) = 2$, then $\gammathin_A(\mu, u) = \mu\ronnn_A$.
\item[(e)]\hypertarget{careful} The family of mappings $\gammathin$ has
the following isometry-equivariance property: for any isometry
$\theta$ of $\R^d$, and for all $(\mu, u) \in\XX\times
[0,1]$,
\[
\theta(\gammathin_A(\mu, u)) =
\gammathin_{\theta(A)}(\theta(\mu), u) .
\]
\end{longlist}
\vspace*{6pt}
\end{proposition}

We shall prove Proposition \ref{specquant} in Section
\ref{quantile}. Property \hyperlink{interest}{(d)} of Proposition
\ref{specquant} will be vital to the proof of Theorem
\ref{result}. It states that whenever $X \ronn_A$ has
exactly one point in its support, $\gammathin_A(X \ronn_A, U)$
will have no points, while whenever $X \ronn_A$ has exactly two
points in its support, $X\ronn_A - \gammathin_A(X\ronn_A, U)$
will have no points. Hence when $X \ronn_A$ has exactly one or
two points the locations of these points provide a possible
source of randomness. The next example will illustrate how
property \hyperlink{interest}{(d)} is exploited and will help to
motivate the proof of Theorem \ref{result}. To make use of
property \hyperlink{interest}{(d)}, we shall need the following
elementary lemma.

Let $\modplus$ denote addition modulo one; that is, for $x,y
\in\R$, let $x \modplus y$ be the unique $z \in[0,1)$ such
that $x+y - z \in\Z$.\vspace*{6pt}
\begin{lemma}[(Adding $\uniforma$ random variables modulo $1$)]
\label{modulo} Let $U_1$ and $U_2$ be $\uniforma$ random
variables that are measurable with respect to the
$\sigma$-fields~$\F_1$ and $\F_2$ and such that $U_1$ is
independent of $\F_2$, and $U_2$ is independent of~$\F_1$. If
$U:= U_1 \modplus U_2$, then $U$ is independent of~$\F_1$, $U$
is independent of $\F_2$ and~$U$ is a $\uniforma$ random
variable.
\end{lemma}
\begin{pf}
The proof\vspace*{-1pt} follows from the Fubini theorem and the fact that for
every $x \in\R$ we have $U_1 \modplus x \eqd U_1$. Let
$E \in\F_2$, and let $Q$ be the joint law of $U_2$ and~$\onee_E$.
Let $B \in\borel$. By symmetry, it is enough to show that $\PP(
\{ U \in B \} \cap E) = \PP(U_1 \in B)\PP(E)$. By the independence of
$U_1$ and $\F_2$, we have
\begin{eqnarray*}
\PP( \{ U \in B \} \cap E)
& = & \int\PP( U_1 \modplus x \in B)i \,dQ(x,i) \\
& = & \int\PP(U_1 \in B)i \,dQ(x,i) \\
& = & \PP(U_1 \in B)\PP(E) .
\end{eqnarray*}
\upqed\end{pf}
\begin{example}[(A monotone map $\phi\dvtx\XX\to\XX$ which maps a Poisson
process~$X$ to another of lower intensity such that $X -
\phi(X)$ is also a Poisson process)] \label{firsteg}
Let $\lx> \ly> 0$. Let $X$ be a Poisson point process on $\R^d$
with intensity $\lx$ and law $\law$. Let ${\mathcal P} =
\{ C_i \}_{i \in\N}$ be an indexed partition of $\R^d$ into
equally-sized cubes, all translates of one another, large enough
so that the Lebesgue measure of each cube is larger than
the constant $K(\lx, \ly)$ from Proposition \ref{specquant}. The
monotone mapping $\phi$ will be defined by specifying
$\phi(\cdot) \ronn_C$ for all $C \in{\mathcal P}$.

Let $U = \{ U_i \}_{i \in\N}$ be a sequence of i.i.d.
$\uniforma$ random variables that are independent of $X$. Let
\[
\Phi(X,U) := \sum_{i \in\N} \gammathin_{C_i}(X, U_i) ,
\]
where $\gammathin$ is the splitting from Proposition
\ref{specquant}. The map $\phi$ will be defined so that
$\phi(X) \eqd\Phi(X, U)$ and $X - \phi(X) \eqd X - \Phi(X,
U)$. By properties \textup{\hyperlink{dist}{(c)}}
and~\textup{\hyperlink{monoquant}{(b)}} of
Proposition \ref{specquant}, we deduce that $\phi(X)$ and
$X - \phi(X)$ are Poisson point processes on $\R^d$ with
intensities $\ly$ and $\lyy$.

If $X(C) = 1$, then we say that $C$ is \textup{one-special}, while
if $X(C) =2$, then we say that $C$ is \textup{two-special}. Let
$k^1$ and $k^2$ be the indices of the one-special and
two-special cubes with the least index, respectively. Note that $\law
$-a.s. $k^1$ and $k^2$ are well defined. Let $Z^1$ be the unique
$X$-point in $C_{k^1}$. Let $Z^2_1$ and~$Z^2_2$ be the two
$X$-points in $C_{k^2}$, where $Z^2_1$ is the one closest to
the origin. Let~$C_\bfz$ be the cube containing the origin.
Fix a measurable function $f_{C_\bfz}\dvtx C_{\bfz} \to[0,1]$ such
that if $V$ is a $\uniformb[C_\bfz]$ random variable, then
$f_{C_\bfz}(V)$ is a $\uniforma$ random variable. For each $C
\in{\mathcal P}$, let $c \in C$ be so that $C-c = C_\bfz$, and
let $f_C\dvtx C \to[0,1]$ be defined via $f_C(x) =
f_{C_\bfz}(x-c)$. Since $X$ is a Poisson point process, it
follows from Lemma~\ref{char} that conditional on $k^1$ we
have that $Z^1$ is a $\uniformb[C_{k^1}]$ random variable.
Moreover it is easy to see that $S^1:=f_{C_{k^1}}(Z^1)$ is in
fact a~$\uniforma$ random variable independent of
\[
\F^1:=
\sigma\bigl(\ind_{[X(C_i) \not= 1]} X \ronn_{C_i}\dvtx i \in
\N\bigr) .
\]

Similarly, it is easy to define $S^2$ as a function of $(Z^2_1,
Z^2_2)$ so that $S^2$ is a~$\uniforma$ random variable
independent of
\[
\F^2:=\sigma\bigl(\ind_{[X(C_i) \not= 2]} X
\ronn_{C_i}\dvtx i \in\N\bigr) ,
\]
namely,
\[
S^2:=
f_{C_{k^2}}(Z^2_1) \modplus f_{C_{k^2}}(Z^2_2) .
\]
To see why
the above definition works, consider the random variables $Y_1$
and $Y_2$, defined as follows. Choose, with a toss of a fair
coin (i.e., independent of $X$), one of $Z^2_1$ or $Z^2_2$ to
be $Y_1$, and let $Y_2$ be so that $\{ Y_1, Y_2 \} =\{ Z^2_1, Z^2_2
\}$. Clearly $Y_1$ and $Y_2$ are independent
$\uniformb[C_{k^2}]$ random variables and $S^2 =
f_{C_{k^2}}(Y_1) \modplus f_{C_{k^2}}(Y_2)$.\vspace*{1pt}

Note that $S^1$ is measurable with respect to $\F_2$, and $S^2$
is measurable with respect to $\F^1$. Let
\[
S:= S^1 \modplus
S^2 .
\]
By Lemma \ref{modulo}, we have that $S$ is
independent of $\F_1$, and $S$ is independent of~$\F_2$. For
all $i \in\N$, let
\[
\phi(X) \ronn_{C_i}:=
\gammathin_{C_i}(X, g_i(S)) ,
\]
where $g_i$ is the
sequence of functions from Lemma \ref{reprod}. By property
\textup{\hyperlink{monoquant}{(b)}} of Proposition~\ref{specquant}, we see that
$\phi$ is monotone. We shall now show that $\phi(X) \eqd
\Phi(X,U)$ and $X - \phi(X) \eqd X- \Phi(X, U)$.

Observe that by property \textup{\hyperlink{interest}{(d)}} of Proposition
\ref{specquant}, for each one-special cube~$C$ we have
\[
\phi(X) \ronn_C = \Phi(X, U) \ronn_C =0 .
\]
Since $S$ is
independent of $\F^1$ and $\{ g_i(S) \}_{i \in\N} \eqd
\{ U_i \}_{i \in\N}$, we have that
\[
\phi(X)=\sum_{ i \in\N} \ind_{[X(C_i) \not= 1]} \gammathin_{C_i}(X,
g_i(S)) \eqd\sum_{i \in\N} \ind_{[X(C_i) \not=1 ]}
\gammathin_{C_i}(X, U_i) = \Phi(X,U) .
\]
Thus $\phi(X) \eqd\Phi(X, U)$.
Similarly, by property \textup{\hyperlink{interest}{(d)}} of Proposition
\ref{specquant}, for each two-special cube $C$ we
have
\[
\bigl(X- \phi(X)\bigr) \big|_C = \bigl(X- \Phi(X, U)\bigr) \big|_C =0 .
\]
Since $S$ is independent of $\F^2$, we have that
\begin{eqnarray*}
X - \phi(X) &=& \sum_{ i
\in\N} \ind_{[X(C_i) \not= 2]}\bigl(X \ronn_{C_i}-
\gammathin_{C_i}(X, g_i(S))\bigr) \\ &\eqd& \sum_{i \in\N}
\ind_{[X(C_i) \not= 2]}\bigl(X \ronn_{C_i} -\gammathin_{C_i}(X,
U_i)\bigr) \\ &=& X - \Phi(X,U).
\end{eqnarray*}
Thus $\phi(X) \eqd\Phi(X, U)$ and $X - \phi(X) \eqd X- \Phi(X,
U)$.
\end{example}

As an aside, one might ask whether the two Poisson processes
$X$ and $X-\phi(X)$ in Example \ref{firsteg} or Theorem
\ref{result} can be made independent of each other, but it
turns out that this is easily ruled out. (It may come as
a~surprise that two dependent Poisson processes can have a sum
that is still a Poisson process; see \cite{jacod}.)
\begin{proposition}
\label{toomuch} There does not exist a monotone map $\phi\dvtx\XX
\to\XX$ such that if $X$ is a homogeneous Poisson point
process on $\R^d$, then $\phi(X)$ and $X - \phi(X)$ are
independent homogeneous Poisson point processes on $\R^d$ with
strictly positive intensities.
\end{proposition}
\begin{pf}
Let $X$ be a Poisson point process on $\R^d$ with intensity
$\lx> 0$. Let $ {\alpha} \in(0,1)$. Toward a contradiction assume
that $\phi(X)$ and $X - \phi(X)$ are independent Poisson point
processes in $\R^d$ with intensities $ {\alpha} \lx$ and $(1-
{\alpha} )\lx$.
Let $\mathfrak R$ and $\mathfrak B$ be independent Poisson
point processes on $\R^d$ with intensities~${\alpha} \lx$ and $(1-
{\alpha} )\lx$. Note that
%
%
\begin{equation}
\label{samedist}
(\mathfrak R, \mathfrak B, \mathfrak R + \mathfrak B) \eqd\bigl(\phi(X),
X- \phi(X), X\bigr) .
\end{equation}
Now let $Z:= \mathfrak R + \mathfrak B$ and let $ B = B(0,1)$, and
consider the events
\[
E := \{ Z(B) = 1 \} \cap\{ \mathfrak R(B) =1 \}
\]
and
\[
E':= \{ X(B) = 1  \} \cap\{ \phi(X)(B) =1 \}.
\]
Clearly,\vspace*{-1pt} $\PP(E \mid Z) = {\alpha} \ind_{[Z(B) =1]}$, but since $E'
\in
\sigma(X)$, we have that $\PP({E'} \mid X) = \onee_{E'}$. Since
$ {\alpha} \in(0,1)$, we conclude that $\PP(E \mid Z ) \not\eqd
\PP({E'} \mid X)$, which contradicts~(\ref{samedist}).
\end{pf}

\subsection*{Outline of the proofs}
Following the lead of Examples \ref{firstfirsteg} and
\ref{firsteg}, we shall introduce an isometry-equivariant
partition of $\R^d$. The partition will consist of globes,
which will be specially chosen balls of a fixed radius,
together with a single unbounded part. The partition will be
chosen as a deterministic function of the Poisson process by a
procedure that does not need to examine the Poisson points
inside the globes. The precise definition of this partition
and its properties are somewhat subtle; see Sections
\ref{selectionrule} and \ref{goodclumping}. The most important
property is that conditional on the partition, the process
restricted to the bounded parts is a Poisson point process that
is independent of the process on the unbounded part. This may
be regarded as an extension of the following property enjoyed
by stopping times for a one-dimensional Poisson process:
Conditional on the stopping time, the process in the future is
a Poisson process independent of the process in the past. The
precise formulation of the property we need may be found in
Proposition \ref{fund}.

To prove Theorem \ref{result}, we shall employ the splitting
from Proposition \ref{specquant} on the bounded parts as in
Example \ref{firsteg}. The Poisson points in the unbounded
part will be split independently of each other with
probabilities $(\frac{\ly}{\lx}, 1 - \frac{\ly}{\lx})$. When
one of the balls of the partition contains exactly one or two
points, the splitting from Proposition \ref{specquant} is
completely deterministic. Thus the locations of these points
provide a source of randomness that can be used to facilitate
the splitting from Proposition \ref{specquant} on the other
balls of the partition, as in Example \ref{firsteg}, and, in
addition, can be used to independently split the points that do
not belong to a bounded part. Of course, we cannot use
randomness precisely as in Example \ref{firsteg} since that
privileges the origin and therefore is not equivariant.
Instead, we use randomness from the available source that is
(essentially) nearest to where it is used.

Aside from some careful bookkeeping to ensure isometry-equivariance,
the two main ingredients for the proof of Theorem
\ref{result} are an isometry-equivariant partition with the
independence property described above and the splitting from
Proposition~\ref{specquant}. Next we focus our discussion on
these two ingredients.

The radius $R$ of the balls of the isometry-equivariant
partition will depend on $(\lx,\ly, d)$. For all $x \in\R^d$
and all $0< s< r$, we define the \textit{shell} centered at~$x$
from $s$ to $r$ to be the set
\[
A(x;s,r): = \{ y \in\R^d\dvtx s \le\|x-y\| \le r\} .
\]
Let $X$ be a Poisson point process on $\R^d$ and $x \in\R^d$. A
single ball of radius~$R$ contained in $B(x, R+10)$ will be
chosen to be a globe (a member of the partition) only if two
properties are satisfied: the shell $A(x; 3R+75 +d,\allowbreak 5R+100 +d)$
contains no $X$-points; and the shell $A(x; R+10, 3R+75+d)$ is
relatively densely filled with $X$-points, that is, every ball of
radius $1/2$ that is contained in $A(x; R+10, 3R+75+d)$ itself
contains an $X$-point. A~minor complication is that the set of
$x\in\R^d$ satisfying these properties is not discrete, but
consists of small well-separated clusters. Each cluster will have
diameter at most $2$ and will be contained in a unique ball of minimum
diameter; the centers of these balls will be the centers of the globes.

The key step in defining the splitting in Proposition \ref{specquant}
is to
construct a coupling of Poisson random variables with the analogous
properties of Proposition \ref{specquant} (save isometry-equivariance). We
shall obtain the joint mass function of the required coupling by
applying a
finite sequence of perturbations to the joint mass function for two
independent Poisson random variables $X$ and~$Y$. Each
perturbation will redistribute the joint probabilities associated with three
consecutive values of each of $X$ and $Y$, while preserving the marginal
distributions of $X$, $Y$ and their sum. See Lemma \ref{split}.

The isometry-equivariant partition used in the proof of Theorem
\ref{result} is used again in the proof of Theorem
\ref{resulttwo}, except that the radius $R$ of the balls will not
depend on $(\lx, \ly, d)$, and we shall set $R=1$; given this
partition, the ideas in Example
\ref{firstfirsteg} can be easily adapted to prove a (weaker)
translation-equivariant version of Theorem~\ref{resulttwo}. It
requires some additional effort to prove Theorem
\ref{resulttwo} in its entirety. The proof of Theorem
\ref{finitary} is not difficult and will follow from the
definitions of the maps in Theorems \ref{result} and
\ref{resulttwo}.

\subsection*{Organization of the paper}
The rest of paper proceeds as follows. In Section
\ref{sectionthick} we prove Theorem \ref{nothick}. This
section is independent of the other sections. Section\vadjust{\eject}
\ref{quantile} is devoted to a proof of Proposition
\ref{specquant}. In Sections \ref{selectionrule} and
\ref{goodclumping} we specify the properties that the
isometry-equivariant partition must satisfy and prove that
such a partition does indeed exist. In Section
\ref{pretheorem} we define some desired properties of a
procedure that assigns randomness from the globes that contain
exactly one or two points to the other globes and to the points of
the unbounded part. The proof of Theorem \ref{result} is given
in Section \ref{def}, and the existence of the procedure that
assigns randomness is proved in Section~\ref{deferproof}. In
Section~\ref{proofoftheoremtwo} we prove Theorem
\ref{resulttwo}. In Section \ref{finitarysection} we prove
Theorem \ref{finitary}. Finally, in Section \ref{openprob} we
state some open problems.

\section{\texorpdfstring{Proof of Theorem \protect\ref{nothick}}{Proof of Theorem 3}}
\label{sectionthick}
\begin{figure}[b]

\includegraphics{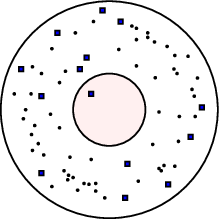}

\caption{The dots are the original points of $X$, and the
squares are points of $\phi(X)\setminus X$. The shaded region
is $B(\bfz,1)$, and the unshaded shell is $A(\bfz;1,M)$. By
selecting subsets of the points in $A(\bfz;1,M)$ uniformly at
random there is nonzero probability that we shall select all
the dots.}
\label{thickpic}
\end{figure}

In this section we shall prove Theorem \ref{nothick}. The proof is by
contradiction. The basic idea is as follows. Let $X$ be a Poisson point
process on $\R^d$ with positive intensity $\lx$ and law $\law$. Let
$\phi\dvtx
\XX\to\XX$ be strongly finitary with respect to $\law$ such that
$\phi(X)$
is a Poisson point process on $\R^d$ with intensity $\ly> \lx$ and
$X \leq
\phi(X)$. Since $\phi(X)$ has greater intensity than $X$, with nonzero
probability we have $X(B(\bfz,1)) = 0$ and $\phi(X)(B(\bfz,1)) \geq
1$. Since
$\phi$ is strongly finitary with respect to $\law$, there is a fixed
deterministic $M$ such that with nonzero probability, we also
have $\phi(X) \ronnn_{B(\bfz,1)} = \phi(X' )\ronnn_{B(\bfz,1)}$,
where $X'$
is equal to $X$ on $B(\bfz, M)$ but is resampled off $B(\bfz,M)$.
Define a
new simple point process $Z$ from $\phi(X)$ by deleting all points in
$B(\bfz, 1)$ and by deleting each point in $[\phi(X) \ronnn_{B(\bfz
,1)^c}]$ independently with probability $\lambda/
\lambda'$
conditional on $\phi(X)$. See Figure \ref{thickpic} for an illustration.
%
%
Since $\phi(X)$ is a Poisson point process, $\phi(X)
\ronnn_{B(\bfz,1)}$ is independent of $\phi(X)
\ronnn_{B(\bfz,1)^c}$, and we may define $Z$ so that it is
independent of $\phi(X) \ronnn_{B(\bfz,1)}$. Since $X \leq
\phi(X)$, there is a nonzero probability that $Z
\ronn_{B(\bfz, M)}=X \ronn_{A(\bfz;1,M)}$. Moreover,
conditional on the event that $X(B(\bfz,1))\,{=}\,0$ and
\mbox{$\phi(X)(B(\bfz,1))\,{\geq}\,1$}, there is a nonzero probability
that $\phi(Z) \ronnn_{B(\bfz,1)} = \phi(X) \ronnn_{B(\bfz,1)}$.
Clearly, this contradicts the independence of $Z$ from $\phi(X)
\ronnn_{B(\bfz,1)}$; the following lemma formalizes this
intuition.
\begin{lemma}
\label{fubini} Let $(S, {\mathcal{S}})$ be a measurable space.
If $X$ and $Y$ are independent random variables taking values in
$S$ and if $A:=\{ y \in S \dvtx\PP(Y=y) > 0 \}$, then $\PP(\{ X = Y\}
\cap\{ Y \in A^c \}) = 0$.
\end{lemma}
\begin{pf}
We apply the Fubini theorem and the independence of $X$ and~$Y$
as follows. Let $\mu_X$ be the law of $X$. Then
\begin{eqnarray*}
\PP(\{ X = Y\} \cap\{ Y \in A^c \}) & = & \PP(\{ X = Y \} \cap\{
X \in A^c \}) \\
& = & \int_{A^c} \PP(Y=x) \,d\mu_X(x) \\
& = & \int_{A^c} 0 \,d\mu_X(x)
= 0 .
\end{eqnarray*}
\upqed\end{pf}

With Lemma \ref{fubini} we can now make the above argument
for Theorem \ref{nothick} precise.
\begin{pf*}{Proof of Theorem \ref{nothick}}
Let $\ly> \lx> 0$. Toward a contradiction, let $X$ be a
Poisson point process on $\R^d$ with intensity $\lx$ and law
$\law$. Let $\phi\dvtx\XX\to\XX$ be a mapping that is strongly
finitary with respect to $\law$ such that $X \leq\phi(X)$ and
$\phi(X)$ is a Poisson point process on $\R^d$ with intensity
$\ly$. Since $X \leq\phi(X)$ and $\phi(X)$ has greater
intensity, we must have that
\[
\PP\bigl( \{ \phi(X)(B(\bfz,1)) \geq1 \} \cap\{ X(B(\bfz,1)) =0
 \} \bigr) > 0 .
\]
Since $\phi$ is strongly finitary, for $\law$-a.e. $\mu\in\XX$,
let $N=N(\mu)$ be the smallest natural
number such that for $\law$-a.e. $\mu' \in\XX$ we have
$\phi(\mu) \ronnn_{B(\bfz,1)} = \phi(\mu')
\ronnn_{B(\bfz,1)}$ whenever $\mu\ronnn_{B(\bfz,N)} =
\mu'\ronnn_{B(\bfz,N)}$. Let
\[
E:= \{ N(X) < M \} \cap
\{ \phi(X)(B(\bfz,1)) \geq1 \} \cap\{ X(B(\bfz,1)) =0 \}
\]
for some $M > 0$. Since $\phi$ is strongly finitary with
respect to $\law$, we have that $\PP(N(X) < \infty) =1$, and
we may choose $M$ so that
%
%
\begin{equation}
\label{greaterthanzero}
\PP(E) > 0 .
\end{equation}
Note that on the event $E$ we have that
\[
\phi(X)
\ronnn_{B(\bfz,1)} = \phi\bigl(X \ronnn_{A(\bfz;1,M) } + W
\ronnn_{B(\bfz,M) ^c} \bigr)\big|_{B(\bfz,1)} ,
\]
where $W$ is
independent of $X$ and has law $P$. Let $U$ be a $\uniforma$
random variable independent of $X$ and $W$. We shall show that
there exists a measurable function
$H\dvtx\XX\times\XX\times[0,1] \to\XX$
such that
%
%
\begin{equation}
\label{contradict}
\PP\bigl( \bigl\{ H\bigl(\phi(X) \ronnn_{A(\bfz;1,M)}, W, U\bigr)
= \phi(X) \ronnn_{B(\bfz,1)} \bigr\} \cap E \bigr) >
0 .
\end{equation}

Define a measurable function $s\dvtx\XX\times[0,1] \to\XX$ such
that if $\mu(\R^d) = \infty$, then $s(\mu, u) = \emptypm$ for
all $u \in[0,1]$, while if $\mu(\R^d) < \infty$, then
$[s(\mu, U)]$ is a~uniformly random subset of
$[\mu]$. Since $X \leq\phi(X)$ and since $U$ is
independent of~$X$, we claim that for any event $E'$ that is measurable
with respect to $X$ and has positive probability,
%
%
\begin{equation}
\label{fubiniverify}
\PP\bigl( \bigl\{s\bigl(\phi(X) \ronnn_{A(\bfz;1,M)}, U \bigr)
= X \ronnn_{A(\bfz;1,M)} \bigr\} \cap E' \bigr) > 0 .
\end{equation}

To verify (\ref{fubiniverify}), let
\[
L:= \int_0 ^1
{\mathbf1}\bigl[s\bigl(\phi(X)\ronnn_{A(\bfz;1,M)}, u\bigr) = X \ronnn_{A(\bfz
;1,M)}\bigr]
\,du .
\]
By the Fubini theorem and the independence of $X$ and $U$, we
have that
\[
\PP\bigl( \bigl\{s\bigl(\phi(X) \ronnn_{A(\bfz;1,M)}, U
\bigr) = X \ronnn_{A(\bfz;1,M)} \bigr\} \cap E' \bigr) = \E L
\onee_{E'}.
\]
Observe that from the definition of $s$
and the fact that $X \leq\phi(X)$, we must have that $L > 0$
$\law$-a.s. Since $\onee_{E'} \geq0$ and $\E\onee_{E'} >
0$, it follows that $\E L \onee_{E'} > 0$.

Hence taking $E' = E$, from (\ref{greaterthanzero}) and
(\ref{fubiniverify}) we have that
%
%
\begin{equation}
\label{contradicttwo}
\PP\bigl(\bigl\{s\bigl(\phi(X) \ronnn_{A(\bfz;1,M)}, U \bigr) = X \ronnn_{A(\bfz
;1,M)} \bigr\} \cap E \bigr) > 0 .
\end{equation}
For all $(\mu,\mu', u) \in\XX\times\XX\times[0,1]$, define
\[
H(\mu,\mu', u):=
\phi\bigl(s\bigl(\mu\ronnn_{A(\bfz;1,M)}, u\bigr)+ \mu' \ronnn_{B(\bfz,M) ^c}
\bigr) \big|_{B(\bfz,1)} .
\]
From (\ref{contradicttwo}),
the definition of $H$ and the definition of $E$, it is obvious
that (\ref{contradict}) holds.

Since $\phi(X)$ is a Poisson point process, $\phi(X)
\ronnn_{B(\bfz,1)}$ and $\phi(X) \ronnn_{A(\bfz;1,M)}$ are
independent, and since $U$ and $W$ are independent of $X$, we have\break
that $\phi(X) \ronnn_{B(\bfz,1)}$ is independent of $H(\phi(X)
\ronnn_{A(\bfz;1,M)}, W,U)$. In addition,\break $\PP(\phi(X)
\ronnn_{B(\bfz,1)} = \mu) = 0$ for all $\mu\in\XX\setminus
\{ \emptypm \}$ and $\phi(X) \ronnn_{B(\bfz,1)} \not= 0$ on the
event~$E$. Thus equation~(\ref{contradict}) contradicts
Lemma \ref{fubini}.
\end{pf*}

\section{\texorpdfstring{Proof of Proposition \protect\ref{specquant}}{Proof of Proposition 8}}
\label{quantile}
The proof of Proposition \ref{specquant} is based on a~spe\-cific coupling of two Poisson random variables.
\begin{lemma}
\label{split} For any $\alpha\in(0,1)$, there exists a $k( {\alpha} )$
such that if $\lambda>k( {\alpha} )$, then there exist random variables
$X$ and $Y$ such that $X$, $Y$ and $X+Y$ have Poisson
distributions with respective means $\alpha\lambda$,
$(1-\alpha)\lambda$ and $\lambda$, and
\[
\PP(Y=0\mid X+Y=1)=1=\PP(X=0\mid X+Y=2) .
\]
\end{lemma}
\begin{pf}
Write $\pi^\gamma_i:=e^{-\gamma} \gamma^i/i!$ for the Poisson
probability mass function. We must find an appropriate joint
mass function for $X$ and $Y$, that is, an element $Q$ of the
vector space $\R^{\N^2}$ with all components nonnegative and
satisfying
%
%
\begin{equation}\label{marg}
\sum_j Q_{i,j}=\pi^{\alpha\lambda}_i ,\qquad
\sum_i Q_{i,j}=\pi^{(1-\alpha)\lambda}_j ,\qquad
\sum_i Q_{i,k-i}=\pi^{\lambda}_k
\end{equation}
and
%
%
\begin{equation}\label{zeros}
Q_{0,1}=Q_{1,1}=Q_{2,0}=0 .
\end{equation}

Let $P\in\R^{\N^2}$ be the mass function for independent
Poisson random variables, that is,
$P_{i,j}:=\pi^{\alpha\lambda}_i\pi^{(1-\alpha)\lambda}_j$, and
note that $P$ satisfies (\ref{marg}) (with $P$ in place of
$Q$). For $s,t\in\N$ define $E^{s,t}\in\R^{\N^2}$ by
$E^{s,t}_{i,j}:=0$ for $(i,j)\notin[s,s+2]\times[t,t+2]$, and
\[
E^{s,t}_{i,j}:=
\begin{tabular}{c|ccc}
$i \diagdown j$ & $t$ & $t+1$ & $t+2$ \\ \hline
$s$ & \hphantom{$-$}0 & $-$1 & \hphantom{$-$}1 \\
$s+1$ & \hphantom{$-$}1 & \hphantom{$-$}0 & $-$1 \\
$s+2$ & $-$1 & \hphantom{$-$}1 & \hphantom{$-$}0 \\
\end{tabular}
\]
and note that $\sum_j E^{s,t}_{i,j}=\sum_i E^{s,t}_{i,j}=\sum_i
E^{s,t}_{i,k-i}=0$.

Now let
\[
Q:=P+P_{0,1}E^{0,0}-(-P_{0,1}+P_{2,0})E^{1,0}-(-P_{0,1}+P_{2,0}+P_{1,1})E^{0,1}
.
\]
From the definition of $Q$, it is easy to verify that
(\ref{zeros}) holds. [The idea is that adding a multiple of
$E^{s,t}$ moves mass from location $(s,t+1)$ to $(s+1,t)$,
without affecting the locations $(i,j)$ with $i+j\leq s+t$.
First we transfer mass~$P_{0,1}$ from location $(0,1)$ to
$(1,0)$; this results in mass $P_{2,0}-P_{0,1}$ at $(2,0)$,
which we then transfer to $(1,1)$; finally we similarly
transfer the current mass at $(1,1)$ to $(0,2)$.] The
equalities in (\ref{marg}) follow from the above observations on
sums involving $P$ and~$E$, so it remains only to check
nonnegativity of $Q$ for $\lambda$ sufficiently large. This
follows by noting that for some $c=c(k,\alpha)>0$ we have
$P_{i',j'}\geq c\lambda P_{i,j}$ whenever $i+j=k$ and
$i'+j'=k+1$; therefore it suffices to take $\lambda$ large
enough compared with
$c(1,\alpha)^{-1},\ldots,c(4,\alpha)^{-1}$.
\end{pf}

For later convenience we next rephrase Lemma \ref{split} in
terms of a mapping that constructs $X$ from $X+Y$.
\begin{corollary}
\label{splitwithf}
For any $\alpha\in(0,1)$, there exists a
$k( {\alpha} )$ such that for $\bar{\lx} > k( {\alpha} )$, there exists
a measurable function $F\dvtx\N\times[0,1] \to\N$ with the
following properties:
\begin{longlist}[(a)]
\item[(a)]
\hypertarget{monof}
For all $(n,u) \in\N\times[0,1]$, we have that $F(n, u) \leq n$.
\item[(b)]
\hypertarget{propinh}
For all $u \in[0,1]$, we have that $F(1,u)=1$ and $F(2,u)=0$.
\item[(c)]\hypertarget{distinh} If $\bar{X}$ is a Poisson random variable
with mean $\bar{\lx}$,
and $U$ is a $\uniforma$ random variable
independent of $\bar{X}$, then $F(\bar{X}, U)$ and $\bar{X} -
F(\bar{X}, U)$ are\vspace*{1pt} Poisson random variables with means
$ {\alpha} \bar{\lx}$ and $(1- {\alpha} )\bar{\lx}$, respectively.
\end{longlist}
\end{corollary}
\begin{pf}
Let $ {\alpha} \in(0,1)$ and $k( {\alpha} )$ be as in Lemma \ref
{split}. Let
$\bar{X}$ be a Poisson random variable with mean $\bar{\lx} >
k( {\alpha} )$, and let $U$ be a $\uniforma$ random variable independent\vadjust{\eject}
of $\bar{X}$. By Lemma \ref{split}, let $X$ and $Y$ be Poisson
random variables with respective means $ {\alpha} \bar{\lx}$ and
$(1- {\alpha} )\bar{\lx}$ such that
$X+Y \eqd\bar{X}$. Define $F$ so that
\[
(\bar{X}, F(\bar{X}, U)) \eqd(X+Y, X) .
\]
\upqed\end{pf}

With Corollary \ref{splitwithf} the proof of Proposition
\ref{specquant} is relatively straightforward, except that
property \hyperlink{careful}{(e)} requires a little care. We next
present some definitions and elementary facts about Poisson
processes that will be useful in the proof and in the rest of
the paper.

Recall that for $\ww\in\XX$, we denote the restriction of
$\ww$ to a set $A \in\borel$ via
\[
\ww\ronn_A(\cdot):=
\ww(\cdot\cap A) .
\]
Recall that \mbox{$\| \cdot\|$} is the Euclidean norm in $\R^d$. We
say that the inter-point distances of a point measure $\mu\in
\XX$ are distinct if for all $x,y,u,v \in[\mu]$ such that
$\{ x,y \} \not= \{ u,v \}$ and $x \not= y$, we have that
$\|x-y\| \not= \|u-v\|$.
\begin{lemma}[(Elementary facts about Poisson point processes)]
\label{interinter}
Let $X$ be a Poisson point process on $\R^d$ with positive intensity
and law $\law$.
\begin{longlist}[(a)]
\item[(a)]
\hypertarget{interc} Let $a \in\R^d$. The distances from the
$X$-points to the point $a$ are distinct $\law$-a.s.
\item[(b)]
\hypertarget{inter}
For all $d \geq1$, the inter-point distances of $X$ are
distinct $\law$-a.s.
\item[(c)]\hypertarget{span} $\law$-a.s., every set of $d$ elements of
$[X]$ has linear span equal to all of~$\R^d$.
\end{longlist}
\end{lemma}
\begin{pf}
The proof follows easily from Lemma \ref{char}.
\end{pf}
\begin{pf*}{Proof of Proposition \ref{specquant}}
Let $X$ be a Poisson point process on $\R^d$ with intensity
$\lx> 0$. Let $ {\alpha} := {\ly}/{\lx}$, and let $k( {\alpha} )$
be defined
as in Corollary~\ref{splitwithf}. Let $K > 0$ be so that
$\bar{\lx}:= K\lx> k( {\alpha} )$. Let $A \in\borel$ have Lebesgue
measure larger than~$K$. Let $\bar{X}:= X(A)$, so that
$\bar{X}$ is a Poisson random variable. Let $F$ be a
function as in Corollary \ref{splitwithf}. Let $U$ be a
$\uniforma$ random variable independent of $X$. Let $g_1, g_2\dvtx
[0,1] \to[0,1]$ be two functions as in Lemma~\ref{reprod} so
that $U_1:=g_1(U)$ and $U_2:=g_2(U)$ are independent $\uniforma$
random variables. Note that by property \hyperlink{monof}{(a)} of
Corollary \ref{splitwithf}, $F(X(A), U_1) \leq X(A)$. We shall
define $\gammathin_A$ so that $[\gammathin_A(X, U)]$ is a
subset of $[X \ronnn_A]$ of size $F(X(A),U_1)$. Moreover,
conditional on $F(X(A), U_1)=j$, each subset of $[X \ronnn_A]$
of size $j$ will be chosen uniformly at random using the randomness
provided by $U_2$. To do this carefully, we shall tag the points
in $[X \ronnn_A]$ and specify a way to use the randomness
provided by $U_2$.

Let $\mu\in\XX$. Consider the following enumeration of the points in
$[\mu\ronnn_A]$. The \textit{center of mass} of a Borel set $C$ with
positive finite Lebesgue measure $\les(C) > 0$ is given by
%
%
\begin{equation}
\frac{1}{\les(C)} \int_C x \,d\les(x) \in\R^d .
\end{equation}
Let $a$ be the center of mass of $A$. We say that $\mu$ admits
the \textit{centric enumeration} on $A$ if $\mu(A) > 0$ and if
the distances from $a$ to the points in $[\mu\ronnn_A]$
are distinct. The centric enumeration on $A$ is given by the
bijection $\iota= \iota_{\mu} \dvtx[\mu\ronnn_A] \to
\{ 1,2,\ldots, \mu(A) \}$, where $\iota(x) < \iota(y)$ iff
$\|x-a\| < \|y-a\|$. Note that by Lemma~\ref{interinter},
part \hyperlink{interc}{(a)}, $X$ admits the centric enumeration on $A$
$\law_{\lx}$-a.s. when $X(A) > 0$.

Now we define some auxiliary functions that, when composed with
$\uniforma$ random variables, yield random variables with certain
distributions. For any set $B$, let $\powerset(B)$ denote the
set of all subsets of $B$. Let $\{ s_{i,j} \}_{j \leq i}$ be a
collection of measurable functions where $s_{i,j} \dvtx[0,1] \to
\powerset(\{ 1,2, \ldots, i \})$ has the
property that if $U'$ is a $\uniforma$
random variable, then $s_{i,j}(U')$ is uniformly distributed
over subsets of size $j$ of $\{ 1,2 ,\ldots,i \}$.

For all $\mu\in\XX$ that do not admit the centric enumeration
on $A$, if $\mu(A) \not= 2$, then set $\gammathin_A(\mu, u) =
\emptypm$ for all $u \in
[0,1]$, and if $\mu(A) =2$, then set $\gammathin_A(\mu, u) = \mu$
for all $u \in
[0,1]$. Otherwise, for $(\mu, u) \in\XX\times[0,1]$, we
proceed as follows. If $\mu(A) = i$, let $\iota\dvtx[\mu\ronnn_A] \to
\{ 1,\ldots, i \}$ be the centric enumeration.
Suppose $F(i,g_1(u))=j$. Define $\gammathin_A(\mu, u)$ to be the
simple point measure with support $\{ x \in[\mu\ronnn_A] \dvtx\iota
(x) \in s_{i,j}(g_2(u)) \}$.

Clearly, by definition, $\gammathin_A$ is monotone and
$\gammathin_A(\mu, u) = \gammathin_A(\mu\ronnn_A, u)$. From
Corollary \ref{splitwithf}, property \hyperlink{distinh}{(c)}, it is
immediate that $\gammathin_A(X,U)(A)$ and $X(A)
-\gammathin_A(X$, $U)(A)$ are Poisson random variables with means
$\ly\les(A)$ and $(\lx- \ly)\les(A)$, respectively. Moreover\vspace*{1pt}
it is easy to check with the help of Lemma~\ref{char} that in
fact $\gammathin_A(X,U)$ and $X\ronn_A - \gammathin_A(X,U)$ are
Poisson point processes on $A$ with intensities $\ly$ and
$\lx-\ly$, respectively. Thus properties \hyperlink{dist}{(c)},
\hyperlink{monotonedist}{(a)} and \hyperlink{monoquant}{(b)} all hold. It is
easy to see that property \hyperlink{interest}{(d)} is also inherited
from property \hyperlink{propinh}{(b)} of Lemma \ref{splitwithf}.
Moreover we have the required property \hyperlink{careful}{(e)} since
we enumerated the points in the support of $\mu\ronnn_A$ in an
isometry-equivariant way via the centric enumeration, while the
functions $g_1, g_2, F, s_{i,j}$ are fixed functions
independent of $\mu$ and $A$.
\end{pf*}

\vspace*{12pt}\section{Selection rules}
\label{selectionrule}

We shall now define an important class of isometry-equivariant
partitions that will have a certain independence property.
Recall that the open \textit{ball} of radius $r$ centered at $x$
is denoted by $B(x,r): = \{ y \in\R^d\dvtx\allowbreak \|x-y\| < r\}$. The
closed ball is denoted by \mbox{$\bar{B}(x, r):= \{ y \in\R^d \dvtx\|x-y\|
\le r\}$}. Let $\goth{F} \subset\borel$ denote the set
of closed subsets of $\R^d$. An \textit{${R}$-selection rule}
is a~mapping $\Psi\dvtx\XX\to\goth{F}$ that has the following
properties:
\begin{longlist}[(a)]
\item[(a)]
\hypertarget{selpp} If $X$ is a Poisson point
process on $\R^d$ with intensity $\lx> 0$ and law~$\law_{\lx}$, then $\law_{\lx}$-a.s. $\Psi(X)$ is a nonempty
union of disjoint closed balls of radius~$R$.
\item[(b)]
\hypertarget{seliso} The map $\Psi$ is isometry-equivariant; that is,
for all isometries $\theta$ of
$\R^d$ and all $\mu\in\XX$, we have that $\Psi(\theta
\ww) = \theta\Psi(\ww)$.
\item[(c)]
\hypertarget{good}
For all $\ww, \ww' \in\XX$, provided $\mu$ and $\mu'$ agree on
the set
%
%
\begin{equation}
\label{hofw}
H(\ww)=H_{\Psi}(\ww):= \biggl( \bigcup_{x \in\Psi(\ww) } \bar
{B}(x,2) \biggr)^c ,
\end{equation}
we have that $\Psi(\ww) = \Psi(\ww')$.
\item[(d)]
\hypertarget{mblty} The map $\Psi$ is measurable; see below
for the precise meaning of this.
\end{longlist}
Let $\Psi$ be an $R$-selection rule, and let $\ww\in\XX$. We
call the connected components of $\Psi(\mu)$ the \textit{globes}
(under $\mu$), and we denote the set of globes by
$\operatorname{Globes}[\Psi(\mu)]$. The \textit{ether} is $\Psi
(\ww) ^c:= \R^d
\setminus\Psi(\mu)$. Note that the set of globes together
with the ether form an isometry-equivariant partition of
$\R^d$.

Note that the set $H(\ww)$ is obtained by first extending
$\Psi(\ww)$ by distance~$2$ and then taking the complement of the
enlarged set.
The idea behind the key condition \hyperlink{good}{(c)} is that $\Psi(\ww)$ is
determined only by the restriction of $\ww$ to~$\Psi(\ww)^c$ [for
technical reasons it is convenient to insist that it is determined even on
the smaller set $H(\ww) \subset\Psi(\mu)^c $, although it seems plausible
that the proof could also be pushed through without this additional
restriction]. This will have the consequence that for a Poisson point
process $X$, conditional on $\Psi(X)$, the process restricted to $\Psi
(X)$ is
still a Poisson point process.
\begin{proposition}
\label{postpone}
For all $d \geq1$ and all $R > 0$, there exists an $R$-selection rule.
\end{proposition}

We postpone the construction of selection rules until Section
\ref{goodclumping}. Sometimes when the value of $R$ is not
important we shall refer to $\Psi$ simply as a selection rule.
The key property of selection rules is the following.
\begin{proposition}[(Key equality)]
\label{fund} Let $X$ and $W$ be independent Poisson point
processes on $\R^d$ with the same intensity. For a selection
rule $\Psi$, the process $Z := W \ronn_{\Psi(X)} + X
\ronn_{\Psi(X) ^c}$ has the same law as $X$ and $\Psi(X) =
\Psi(Z)$.
\end{proposition}

Proposition \ref{fund} states that conditional on $\Psi(X)$,
not only is $X\ronnn_{\Psi(X)}$ a~Poisson point process on
$\Psi(X)$, it is also independent of $X\ronnn_{\Psi(X)^c}$.

\subsection*{Some remarks on measurability}
It will be obvious from our construction of selection rules
that measurability will not be an issue. However, for the sake
of completeness and since we want to prove Proposition
\ref{fund} before providing the explicit construction of
selection rules, we assign the \textit{Effros} $\sigma$-algebra
to $\goth{F}$. For each compact set $K \in\borel$, let
$\goth{F}_K:=\{ F \in\goth{F} \dvtx F \cap K \not= \varnothing \}$.
The Effros $\sigma$-algebra for $\goth{F}$ is generated by the
sets $\goth{F}_K$ for all compact sets $K \in\borel$. Let
$(\X, \F, \PP)$ be a probability space. We call a measurable
function $\mathcal{X}\dvtx\X\to\goth{F}$ a \textit{random closed
set}. Thus if $X$ is a Poisson point process and $\Psi$ is a
selection rule, then $\Psi(X)$ is a random closed set. We
shall not need to use any results from the theory of random
closed sets; we refer the interested reader to
\cite{MR2132405} for background.

\subsection*{\texorpdfstring{Remarks on the proof of Proposition \protect\ref{fund}}{Remarks on the proof of Proposition 16}}
It is immediate from property~\hyperlink{good}{(c)} that $\Psi(X) =
\Psi(Z)$. The isometry-equivariance of selection rules
[property \hyperlink{seliso}{(b)}] will not play a role in the proof
of Proposition \ref{fund}. For the purposes of the following
discussion, let us assume that $\Psi$ does not have to satisfy
property \hyperlink{seliso}{(b)}. Temporarily suppose instead that
$\Psi$ satisfies the following additional requirement:
\begin{longlist}[(b$'$)]
\item[(b$'$)] There exists a \textit{fixed} Borel set $D$
such that if $X$ is a Poisson point process on $\R^d$,
then $\Psi(X) \subset D \subset H_\Psi(X) ^c$.
\end{longlist}
For any random variable $Y$, we let $\sigma(Y)$ be the $\sigma
$-algebra generated by $Y$. By property \hyperlink{good}{(c)} in the definition
of a selection rule,
it is easy to see that~$\Psi(X)$ is~$\sigma(X\ronnn_{D^c})$-measurable. Since $X\ronnn_{D}$ and
$X\ronnn_{D^c}$ are independent, we have that
\[
W \ronn_{D \cap
\Psi(X)} + X \ronn_{D \cap\Psi(X) ^c} \eqd X\ronn_{D} ,
\]
where $W \eqd X$ and $W$ is independent of $X$. Moreover, one
can verify (see Lem\-ma~\ref{prefund} below) that
%
%
\begin{equation}
\label{fundverify}
W \ronn_{D \cap\Psi(X)} + X \ronn_{D \cap\Psi(X) ^c} + X\ronn
_{D^c} \eqd X.
\end{equation}
If $\Psi$ satisfies condition (b$'$), then the left-hand side
of (\ref{fundverify}) equals $W \ronn_{\Psi(X)} + X
\ronn_{\Psi(X) ^c}$, and Proposition \ref{fund} follows.

The above argument suggests that to prove Proposition
\ref{fund}, we should examine events where $\Psi(X)$ is
contained in some deterministic set. However, in general, such
events will have probability zero. We can overcome this
problem by considering events where for some bounded Borel set
$A$, we have that $\Psi(X) \cap A$ is contained in some
deterministic set. For each bounded Borel set $A$, Lemma
\ref{grid} below specifies some additional useful properties
that we require of such events.
\begin{lemma}
\label{grid} Let $X$ be a Poisson point process on $\R^d$ with
positive intensity. Let $\Psi$ be an $R$-selection rule, and
let $H$ be defined as in (\ref{hofw}). Let $A$ be a bounded
Borel set. There exists a finite set $F$, a collection of
disjoint events $\{ E( {\alpha} ) \}_{ {\alpha} \in F}$ and a
collection of
bounded Borel sets $\{ D( {\alpha} ) \}_{ {\alpha} \in F}$ with the
following
properties:
\begin{longlist}
\item
\hypertarget{onE}
For all $ {\alpha} \in F$, on the event $E( {\alpha} )$,
we have
that
\[
\Psi(X) \cap A \subset D( {\alpha} ) \subset H(X)^c .
\]
\item
\hypertarget{Emeas} For all $ {\alpha} \in F$, the event $E( {\alpha
} )$ is
$\sigma(X\ronn_{D( {\alpha} )^c })$-measurable.
\item
\hypertarget{probone} The disjoint union $\bigcup_{ {\alpha} \in F}
E( {\alpha} )$ is an event of probability one.
\end{longlist}
\end{lemma}

We shall prove this later.

The following lemma will be useful in the proof of Proposition
\ref{fund}. In particular, it justifies equation
(\ref{fundverify}) when $\Psi$ satisfies condition (b$'$). The
lemma is a technical generalization of the fact if $X$ and $W$
are two independent Poisson point processes on $\R^d$ with the
same intensity, then for all $s \in\borel$ we have
%
%
\begin{equation}
\label{fixeds}
W \ronnn_s + X \ronnn_{s^c} \eqd X.
\end{equation}
\begin{lemma}
\label{prefund} Let $\gothh{X}$ and $\gothh{W}$ be independent
homogeneous Poisson point processes with equal intensity on some Borel set
$\gothh{D} \subset\R^d$. Let $\gothh{T}$ be a random closed set,
and let
$\gothh{S}:=\gothh{T} \cap\gothh{D}$. Let $\gothh{Y}$ be any point process
and let $\gothh{V}$ be an event. Let $\gothh{S'} := \gothh{D}
\setminus
\gothh{S}$. If $(\gothh{X},\gothh{W})$ is independent of
$(\gothh{S},\gothh{Y},\gothh{V})$, then for all measurable sets of point
measures $\A\in\mathcal{M}$,
\[
\PP( \{ \gothh{X} + \gothh{Y} \in\A \} \cap{\gothh{V}}
) =\PP(
\{ \gothh{W} \ronn_{\gothh{S}} + \gothh{X} \ronn_{\gothh{S'}}+
\gothh{Y} \in\A \} \cap{\gothh{V}} ) .
\]
\end{lemma}
\begin{pf}
Let $\mu_{\gothh{X}}$ be the law of $\gothh{X}$, and let $Q$ be
the joint law of $\gothh{S},\gothh{Y}$ and~$\onee_{\gothh{V}}$.
From (\ref{fixeds}) it is easy to see that for all Borel $s
\subset\gothh{D}$, and for all $\A' \in\mathcal{M}$,
%
%
\begin{eqnarray}
\label{fixedss}
\int\ind_{x \ronnn_s + x \ronnn_{s'} \in\A'} \,d\mu_{ \gothh{X}}(x)
&=& \PP(\gothh{X} \in\A') \nonumber\\
&=& \PP( \gothh{W} \ronnn_s + \gothh{X} \ronnn_{s'} \in\A')
\\
&=& \iint\ind_{[w \ronnn_s + x \ronnn_{s'} \in\A']}
\,d\mu
_{\gothh{X}}(w)
\,d\mu_{\gothh{X}}(x).\nonumber
\end{eqnarray}
Let $\A\in\mathcal{M}$ and $ L:=\PP( \{ \gothh{X} + \gothh
{Y} \in\A\} \cap{\gothh{V}} )$. By the independence
of $\gothh{X}$ and $(\gothh{S},\gothh{Y},\gothh{V})$, we have
that
%
%
\begin{eqnarray}
\label{fixedswithy}
L &=& \iint\ind_{[x + y \in\A]} v \,d\mu_{\gothh{X}}(x)
\,dQ(s,y,v) \nonumber\\[-8pt]\\[-8pt]
&=& \int\biggl(\int\ind_{[x \ronnn_s + x \ronnn_{s'} + y \in\A]}
\,d\mu_{\gothh{X}}(x)\biggr) v \,dQ(s,y,v).\nonumber
\end{eqnarray}
Applying (\ref{fixedss}) to (\ref{fixedswithy}), we obtain that
%
%
\begin{equation}
\label{fullint}
L = \iiint\ind_{[w \ronnn_s +x \ronnn_{s'} + y \in
\A]} v
\,d\mu_{\gothh{X}}(w)\,d\mu_{\gothh{{X}}}(x) \,dQ(s,y,v) .
\end{equation}
Since $\gothh{X}$ and $\gothh{W}$ are independent, and
$(\gothh{X},\gothh{W})$ and $(\gothh{S},\gothh{Y},\gothh{V})$
are independent, we easily recognize that the right-hand side
of equation (\ref{fullint}) is equal to $\PP( \{ \gothh{W}
\ronn_{\gothh{S}} +\gothh{X} \ronn_{\gothh{S'}} + \gothh{Y} \in
\A \} \cap{\gothh{V}} )$.
\end{pf}

With the help of Lemmas \ref{grid} and \ref{prefund} we now
prove Proposition \ref{fund}.
\begin{pf*}{Proof of Proposition \ref{fund}}
Let $X,W\dvtx\X\to\XX$ be independent Poisson point processes on
$\R^d$ with the same intensity, defined on the probability
space $(\X, \F, \PP)$. We shall use $\omega$ to denote an
element of the probability space, and during this proof
$X(\omega)$ will denote the point measure that is the image of
$\omega$ under the random variable $X$ (not ``the number of
$X$-points in $\omega$''). Let $\Psi$ be an $R$-selection rule,
and let $Z:= W \ronn_{\Psi(X)} + X \ronn_{\Psi(X) ^c}$. Let $\A
\in\mathcal{M}$. It suffices to show that $\PP(X \ronnn_A \in
\mathcal{A}) = \PP(Z \ronnn_A \in\mathcal{A})$ for all
bounded Borel sets~$A$. Let $A$ be a bounded Borel set, and let
$\{ E( {\alpha} ) \}_{ {\alpha} \in F}$ and $\{ D( {\alpha} ) \}_{
{\alpha} \in F}$ be
collections of events and subsets of $\R^d$ that satisfy the
conditions of Lemma~\ref{grid}. We shall show that for all
$\alpha\in F$,
%
%
\begin{equation}
\label{ealpha}
\PP\bigl( {\{ X \ronnn_A \in\mathcal{A} \} \cap E(\alpha) } \bigr)
= \PP\bigl( {\{ Z \ronnn_A \in\mathcal{A}  \} \cap E(\alpha)} \bigr) .
\end{equation}
By summing over all $\alpha\in F$, we can then conclude by
property \hyperlink{probone}{(iii)} of Lem\-ma~\ref{grid} that $\PP(X \ronnn_A
\in\mathcal{A}) = \PP(Z \ronnn_A \in\mathcal{A})$.
Let us fix $\alpha\in F$, and set $E := E(\alpha)$ and $D
:=D(\alpha)$. Observe that for all $\w_1, \w_2 \in E$, we
have $\Psi(X(\w_1)) = \Psi(X(\w_2))$ whenever $X(\w_1) =
X(\w_2)$ on $D^c$. This follows from property \hyperlink{good}{(c)}
in the definition of a selection rule and property \hyperlink{onE}{(i)}
of Lemma \ref{grid}. Now define $S := \Psi(X \ronnn_{D^c})$.
Clearly, $S$ is $ \sigma( X \ronnn_{D^c})$-measurable, and on
the event~$E$, we have that $S = \Psi(X)$. Since $X$ is a
Poisson point process, we have that~$X\ronnn_{D \cap A}$ is
independent of $ X\ronnn_{ D^c \cap A}$. Also, by property
\hyperlink{Emeas}{(ii)} of Lemma \ref{grid} we have that $E \in
\sigma( X \ronnn_{D^c})$. See Figure \ref{nong}
for an illustration.

%
%
\begin{figure}[b]

\includegraphics{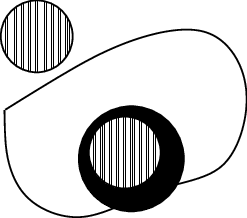}

\caption{An illustration of the deterministic sets
$A$, $D$ and the random set $S$, on the event~$E$ which depends only on
$X \ronnn_{D^c}$. The set $A$ is the large enclosed space, $D$ is the
black disc and $S$ is the union of the hatched discs. The hatched disc
contained in $D$ is $S \cap A$ and its location within $D$ depends only
on $X \ronnn_{D^c}$.} \label{nong}
\end{figure}

By applying Lemma \ref{prefund}
with the following substitutions:
\begin{eqnarray*}
\gothh{D} &=& D \cap A , \qquad\gothh{X}=X \ronnn_{D \cap A} ,\qquad
\gothh{W} = W \ronnn_{D \cap A} , \\
\gothh{T} &=& S ,\qquad
\gothh{Y} = X \ronnn_{D^c \cap A} ,\qquad
\gothh{V} = \onee_ {E} ,
\end{eqnarray*}
it is easy to check that
\[
\PP( \{ X\ronnn_A \in\A \} \cap E
) = \PP( \{ W \ronnn_{D \cap A \cap S} + X \ronnn_{D \cap
A \cap S^c} + X\ronnn_{D ^c \cap A} \in\A \} \cap E
) .
\]
Thus from the definition of $S$ and property
\hyperlink{onE}{(i)} of Lemma \ref{grid}, we have that
\[
\PP( \{ X\ronnn_A \in\A \} \cap E ) = \PP\bigl( \bigl\{ W
\ronnn_{\Psi(X) \cap A} + X \ronnn_{ \Psi(X)^c \cap A} \in\A \bigr\}
\cap E \bigr) .
\]
By the
definition of $Z$, we see that we have verified equation
(\ref{ealpha}) as required.
\end{pf*}

It remains to prove Lemma \ref{grid}.
\begin{pf*}{Proof of Lemma \ref{grid}}
We need some preliminary definitions. The open \textit{cube} of
side length $2r$ centered at the origin is the set $(-r,r)^d$.
The \textit{diameter} of a set $A \subset\R^d$ is ${\sup_{x,y
\in A} }\|x-y\|$. Let $X$ be a Poisson point process on~$\R^d$,
and let $\Psi$ be an $R$-selection rule. Recall that by
property \hyperlink{selpp}{(a)} in the definition of a selection rule,
all globes are balls of radius $R$. Fix a bounded Borel set~$A$.
Let $A' := \bigcup_{x \in A} B(x,2R)$. Let $\{ c_i \}_1
^N$ be a collection of disjoint cubes of diameter $\frac{1}{2}$
such that their union contains the set $A'$. Thus, some cubes
may not be open. Let $a_i \in c_i$ be the centers of the
cubes. Let $F_i$ be the event that the center of some globe
(under $X$) is an element of the cube~$c_i$. For a~binary
sequence $\alpha\in\{ 0,1 \} ^N$ of length $N$, define
\[
E(\alpha):= \biggl(\mathop{\bigcap_{1 \leq i \leq N:}}_{\alpha(i)
=1} F_i\biggr) \cap
\biggl( \mathop{\bigcap_{1 \leq i \leq N:}}_{\alpha(i) =0} F_i ^c
\biggr) .
\]
Set $F:=\{ {\alpha} \in\{ {0,1} \}^{N} \dvtx\PP(E( {\alpha} )) > 0 \}$.
Note that
the events $\{ E(\alpha) \}_{\alpha\in F}$ are disjoint and
their union over all $\alpha$ is an event of probability $1$,
so that condition \hyperlink{probone}{(iii)} is satisfied. Note that if
$x \in\R^d$ and $\|x-a_i\| \leq\frac{1}{2}$, then
%
%
\begin{equation}
\label{easytriangle}
\bar{B}(x,R) \subset B(a_i, R+1) \subset\bar{B}(x,R+2) .
\end{equation}
Define
\[
D(\alpha):= \mathop{\bigcup_{1 \leq i \leq N:}}_
{\alpha(i) =1} B(a_i, R+1) .
\]
Since every globe that intersects $A$ has its center lying at
distance at most~$R$ from~$A$, every globe that intersects $A$
must have a center in some cube~$c_i$. By definition, for
every $ {\alpha} \in F$, on the event $E( {\alpha} )$ we see from~(\ref{easytriangle}) and~(\ref{hofw}) that
\[
\Psi(X) \cap A
\subset D( {\alpha} ) \subset H(X)^c ,
\]
since the diameter of each
cube $c_i$ is $\frac{1}{2}$. See Figure \ref{globeapprox}
for an illustration.

%
%
\begin{figure}

\includegraphics{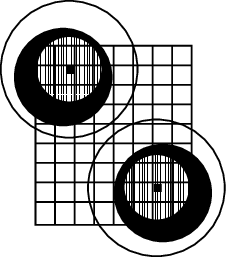}

\caption{The grid is an illustration of the cubes $c_i$. The
black squares are the centers of the globes. The hatched discs
are the globes, the union of the black discs is the set
$D( {\alpha} )$ that contains the globes intersecting $A$ and the
area contained in the largest circles is part of the set
$H(X)^c$.}
\label{globeapprox}
\end{figure}

Thus condition \hyperlink{onE}{(i)} is satisfied. Observe that for
each $ {\alpha} \in F$, we have that $E( {\alpha} ) \in\sigma(X
\ronnn_{D( {\alpha} )
^c})$ by property \hyperlink{good}{(c)} in the definition of a selection
rule, so that condition \hyperlink{Emeas}{(ii)} is also satisfied.
\end{pf*}

Proposition \ref{fund} will be instrumental in proving
Theorems \ref{result} and \ref{resulttwo}. In
Corollary~\ref{fresh} below, we make an important step in this
direction
by constructing a splitting that involves different mechanisms\vadjust{\eject}
on the globes and on the ether. Before stating this result, we
need some preliminary definitions. In particular, recall the
elementary fact that if each point of a Poisson point process
$X$ with intensity $\lx$ is deleted independently of all others
with probability $\frac{\ly}{\lx}$, where $\ly<\lx$, then the
remaining points and the deleted points form independent
Poisson point processes with intensities $\ly$ and $\lx- \ly$.
To facilitate later variations on this theme, we shall give a
very explicit version of this fact.

Sometimes it will be convenient to specify a well ordering of
the sets $[\mu]$, $[\mu\ronnn_{\Psi(\mu) ^c}]$ and
$\operatorname{Globes}[\Psi(\mu)]$. This can be done in the
following way.
Consider the ordering $\prec$ on $\R^d$ in which $x \prec y$
iff $\|x\| < \|y\|$ or iff $\|x\| = \|y\|$ and $x$ is less than
$y$ in the lexicographic ordering of $\R^d$. Thus we can well
order $[\mu]$ and $[\mu\ronnn_{\Psi(\mu) ^c}]$ via
$\prec$ and well order $\operatorname{Globes}[\Psi(\mu)]$ by well
ordering the
centers of the globes via $\prec$. We shall call $\prec$ the
\textit{radial ordering}.

Define $\relstandthin= \relstandthin_{(\lx, \ly)}\dvtx\R^d \times
[0,1] \to\XX$ via
%
%
\begin{equation}
\label{coin}
\relstandthin(x,u):= \ind_{[u \leq{\ly/\lx} ]}\delta_x .
\end{equation}
Define $\standthin= \standthin_{(\lx, \ly)}\dvtx\XX\times[0,1]
\to\XX$ by
%
%
\begin{eqnarray}
\label{standthin}
\standthin_{(\lx, \ly)}(\mu, u):=
\sum_{i=1} ^{\infty} \relstandthin_{(\lx, \ly)}(x_i, g_{i}(u)),\qquad
(\mu,u)\in\XX\times[0,1] ,
\end{eqnarray}
where $\{ x_i \}_{i=1} ^{\infty}$ is $[\mu]$ ordered by
$\prec$ and the $g_i$ are from Lemma \ref{reprod}. We shall
call $\standthin$ the \textit{standard splitting}. The following
fact is elementary.
\begin{lemma}[(Independent splitting)]
\label{standoverset} If $X$ is a Poisson point process on~$\R^d$
with intensity $\lx$, and $U$ is a $\uniforma$ random
variable independent of $X$, then for all $A \in\borel$ and
for all $\ly< \lx$, we have that $\standthin_{(\lx,
\ly)}(X\ronnn_A, U)$ and $X\ronnn_A - \standthin_{(\lx,
\ly)}(X\ronnn_A, U)$ are independent Poisson point processes on
$A$ with intensities $\ly$ and $\lx- \ly$, respectively.
\end{lemma}
\begin{corollary}
\label{fresh} Let $X$ be a Poisson point process on $\R^d$ with
intensity~$\lx$, and let $\ly< \lx$. Let $\gammathin$ be the
splitting from Proposition \ref{specquant}. Let $\Psi$ be an
$R$-selection rule, where the Lebesgue measure of $B(\bfz, R)$
is larger than that of the constant $K(\lx, \ly)$ from
Proposition \ref{specquant}. Let $\{ b_i \}_{i \in\Z^{+}} =
\operatorname{Globes}[\Psi(X)]$, where we have ordered the globes via
the radial
ordering. Let $U$ be a $\uniforma$ random variable
independent of $X$, and let $g_i\dvtx[0,1] \to[0,1]$ be a sequence
of functions as in Lemma \ref{reprod}. The mapping
$\Phi=\Phi_{(\lx, \ly)}$ defined by
\[
\Phi(X, U) := \sum_{i \in\Z^{+}}\gammathin_{(b_i,\lx,\ly
)}(X\ronn_{b_i},g_i(U)) + \standthin_{(\lx,\ly)}\bigl(X \ronn_{\Psi(X)
^c}, g_0(U)\bigr)
\]
is a splitting such that $\Phi(X, U)$ and $X - \Phi(X, U)$ are
Poisson point processes with intensities $\ly$ and $\lx-\ly$, respectively.
\end{corollary}
\begin{pf}
The inequality $\Phi(X, U) \leq X$ is obvious from the
definition of~$\Phi$, so we just need to check that $\Phi(X,U)$
and $X - \Phi(X, U)$ have the right distributions. This is
made possible via Proposition \ref{fund}. Let $W$ be a
Poisson point process on $\R^d$ with intensity $\lx$ that is
independent of $X$ and $U$. Let~$U_1$, $U_2$ be independent
$\uniforma$ random variables that are also independent of $X$
and $W$. From the definition of~$\standthin$, it is easy to
see that
\[
\standthin_{(\lx, \ly)}\bigl(W \ronnn_{\Psi(X)} + X \ronnn_{\Psi(X)
^c}, U_1\bigr)
\eqd\standthin_{(\lx, \ly)}\bigl(W \ronn_{\Psi(X)}, U_1\bigr) + \standthin
_{(\lx, \ly)}\bigl(X \ronn_{\Psi(X) ^c}, U_2\bigr) ,
\]
since the ordering of the points of $W \ronnn_{\Psi(X)} + X
\ronnn_{\Psi(X) ^c}$ is irrelevant as long as the ordering is
independent of $U_1$ and $U_2$. By Proposition \ref{fund},
we have that $X \eqd W \ronnn_{\Phi(X)} + X \ronnn_{\Phi(X)
^c}$, so we obtain that
%
%
\begin{equation}
\label{stand}
\standthin_{(\lx, \ly)}(X,U_1) \eqd
\standthin_{(\lx, \ly)}\bigl(W \ronn_{\Psi(X)}, U_1\bigr)
+ \standthin_{(\lx, \ly)}\bigl(X \ronn_{\Psi(X) ^c}, U_2\bigr) .
\end{equation}
From property \hyperlink{dist}{(c)} of Proposition \ref{specquant} and
Lemma \ref{standoverset}, it is easy to see that for any $A
\in\borel$ with finite Lebesgue measure larger than $K$, we have
\[
\gammathin_A (W \ronnn_{A}, U_1) \eqd\standthin(W \ronnn_{A},
U_1) .
\]
Moreover, since $X$ and $W$ are independent, it follows that
%
%
\begin{equation}
\label{properthin}
\PP\biggl(
\sum_{ i \in\Z^{+}} \gammathin_{b_i} (W\ronnn_{b_i}, g_i(U_1)) \in
\cdot\Bigm| X \biggr)
= \PP\bigl( \standthin\bigl(W \ronnn_{\Psi(X)}, U_1\bigr) \in\cdot\mid
X\bigr) .
\end{equation}
(Recall that $\{ g_i(U) \}_{i \in\N}$ is a sequence of i.i.d.
$\uniforma$ random variables.) Clearly by Proposition
\ref{fund} and the definition of $\Phi$, we have
\begin{eqnarray*}
\Phi(X, U) &\eqd& \Phi\bigl( W\ronnn_{\Psi(X)} + X\ronnn_{\Psi(X)^c},
U \bigr) \\
&\eqd& 
\sum_{ i \in\Z^{+}} \gammathin_{b_i} (W\ronnn_{b_i}, g_i(U_1))
+ \standthin\bigl(X\ronnn_{\Psi(X) ^c}, U_2\bigr) .
\end{eqnarray*}

From equation (\ref{properthin}) and the fact that $X$ and
$W$ are independent, it is easy to verify that
%
%
\begin{equation}
\label{standd}
\Phi(X, U) \eqd\standthin\bigl(W \ronn_{\Psi(X)}, U_1\bigr)
+ \standthin\bigl(X\ronnn_{\Psi(X) ^c}, U_2\bigr) .
\end{equation}
Putting\vspace*{1pt} (\ref{stand}) and (\ref{standd}) together, we obtain
that $\Phi(X, U) \eqd\standthin_{(\lx, \ly)}(X, U_1)$. Thus
from Lemma \ref{standoverset} we have verified that $\Phi(X,
U)$ is a Poisson point process of intensity $\ly$.

The proof that $X - \Phi(X, U)$ is a Poisson point process of
intensity \mbox{$\lx- \ly$} follows by the same argument since
$\gammathin$ is a splitting by Proposition \ref{specquant}
and~$\standthin$ is a splitting by Lemma \ref{standoverset}.
\end{pf}

Let us remark that for Corollary \ref{fresh}, in order for
$\Phi$ to be a splitting we must apply the splitting
$\gammathin$ in all the globes and not just the globes that contain
exactly one or two points.
For example, if $X$ is a Poisson point process on a bounded
Borel set $B$, the following procedure will not result in a
splitting: apply $\gammathin$ if there are exactly one or two
$X$ points, otherwise apply $\standthin$.

Before we begin the proof of Theorem \ref{result}, we first
provide a construction of selection rules along with some other
minor constructions that will be needed.
%
\section{Construction of selection rules}
\label{goodclumping}
\begin{figure}

\includegraphics{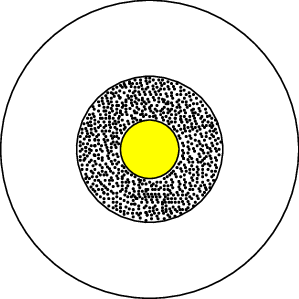}

\caption{An illustration of a pre-seed. The outer shell (the
empty shell) contains no $X$-points. The intermediate shell
(the halo) is relatively densely filled with $X$-points. The
shaded area is unspecified in terms of $X$.}
\label{globe}
\end{figure}

Fix $d \geq1$ and $R > 0$. We shall now construct an
$R$-selection rule. We need some preliminary definitions.
Recall the definition of the \textit{shell},
\[
A(x;s,r): = \{ y \in\R^d\dvtx s \le\|x-y\| \le r \}.
\]
Let $X$
be a Poisson point process on $\R^d$ with intensity $\lx> 0$
and law $\law_{\lx}$. A~point $x \in\R^d$ is called a
\textit{pre-seed} if $B(x,5R + 100\plusd)$ has the following two
properties:
\begin{longlist}[(a)]
\item[(a)]
\hypertarget{empty}
$X(A(x;3R+75\plusd, 5R+100\plusd))=0$;
\item[(b)]
\hypertarget{halo}
for every open ball $B$ of radius
$\frac{1}{2}$ satisfying $B \subset A(x; R+10,
3R+75\plusd)$,
we have $X(B) \geq1$.
\end{longlist}
Given $\mu\in\XX$, we also say that $x$ is pre-seed under
$\mu$ if \hyperlink{empty}{(a)} and \hyperlink{halo}{(b)} hold with~$X$
replaced by $\mu$. If $x$ is a pre-seed, we call $A(x;
3R+75\plusd, 5R+100\plusd)$ the associated \textit{empty shell} and
$ A(x; R+10, 3R+75\plusd)$ the associated \textit{halo}. Clearly
pre-seeds exist $\law_{\lx}$-a.s. An $R$-selection rule will be
defined so that its globes will be balls of radius $R$
contained in $B(x, R+10)$ for some pre-seed~$x$. See Figure
\ref{globe} for an illustration of a pre-seed.

%

Observe that if $x,y \in\R^d$ are pre-seeds, then $\|x-y\|
\notin(2, 2(3R+63)\plusd))$; otherwise the empty shell of one
pre-seed would intersect the halo of the other in such a way as to
contradict the definition of a pre-seed. Also note that the width of
the empty shell is chosen to be greater than $2(R +10)$; this is needed
in the special case $d=1$ to ensure that if $x,y \in\R^d$ are
pre-seeds, then $\|x-y\|
\notin(2, 2(3R+63)\plusd))$. We say that two
pre-seeds $x,y$ are \textit{related} if $\|x-y\| \le2$. This
gives an equivalence relation on the pre-seeds.

We next associate with each
equivalence class a single point in an isometry-equivariant way. Let
$C$ be
an equivalence class of pre-seeds under $\mu$. Observe that $C$ is
contained in some ball of radius $2$, which also contains a unique
point $c \in\R^d$ that is the center of the ball with the smallest
radius that contains $C$; we declare that $c$ is a \textit{seed}. Note
that $c$ might not be a pre-seed (but it has
properties similar to a pre-seed).

If $c$ is a seed (under $\mu$), we call $\bar{B}(c, R)$ a
\textit{globe} (under $\mu$). Define the mapping $\Psi_{R}\dvtx\XX
\to\goth{F}$ by stipulating that for each $\mu\in\XX$, the
Borel set $\Psi_{R}(\ww)$ is the union of the set of globes
under $\mu$.
Given any two
seeds, it is easy to see that their globes do not intersect.
Thus the definition of a globe given here is consistent with
the definition of a globe given in Section
\ref{selectionrule}.

Next, we show that for $R > 0$, the mapping $\Psi_{R}$ is a
selection rule, thus proving Proposition \ref{postpone}.
\begin{lemma}
\label{samepreseed} Let $\Psi= \Psi_{R}$ be the mapping
defined above. For all $\ww, \ww' \in\XX$, if $\ww= \ww'$ on
$H(\ww):={ (\bigcup_{x \in\Psi(\ww) } \bar{B}(x,2)
) ^c }$, then $\ww$ and $\ww'$ have the same pre-seeds.
\end{lemma}
\begin{pf}
Assume that $\ww$ and $\ww'$ agree on $H(\ww)$. Let $z \in
\R^d$ be a pre-seed under~$\ww$. We claim that
\[
\ww
\ronn_{A(z;R+10,5R+100\plusd)} =
\ww'\ronn_{A(z;R+10,5R+100\plusd)}
\]
from which we deduce
that $z$ is also a pre-seed under $\ww'$.

Let $C(z)$ be the equivalence class of pre-seeds to which $z$
belongs, and let~$c$ be the corresponding seed.
Since $c$ has distance at most $4$
from $z$, and $z$ has distance at least $2(3R+63\plusd)$ from any
pre-seed (under $\ww$) not in $C(z)$, we have that $c$ has
distance at least $2(3R + 63 \plusd) -4$ from any pre-seed (under~$\ww$)
not in $C(z)$. Let $m >0$ be the minimal distance from
$c$ to another seed (under $\ww$). Clearly $m \geq
2(3R+63\plusd) - 8$. Since $\ww= \ww'$ on $H(\ww)$, we have
that $\ww\ronn_{A(c;R+2,m-R-2)} = \ww'
\ronn_{A(c;R+2,m-R-2)}$. Since $z$ has distance at most $4$
from $c$, clearly $\ww\ronn_{A(z;R+10,5R+100\plusd)} = \ww'
\ronn_{A(z;R+10,5R+100\plusd)}$, as required.
\end{pf}
\begin{pf*}{Proof of Proposition \ref{postpone}}
Let $R > 0$ and $d \geq1$. We shall now check that $\Psi= \Psi_{R}$
is indeed an $R$-selection rule.

Property \hyperlink{selpp}{(a)}: Let
$\law$ be the law of a Poisson point process on $\R^d$ with
positive intensity. Note that pre-seeds occur $\law$-a.s.
Therefore, we have that seeds occur $\law$-a.s. and
$\Psi(\mu) \not= \varnothing$ for $\law$-a.e. $\mu$. Also, by
definition, if $\Psi(\mu) \not= \varnothing$, then $\Psi(\mu)$ is a
disjoint union of balls of radius $R$.

Property \hyperlink{seliso}{(b)}: Let $\theta$ be an isometry of
$\R^d$. If $x \in\R^d$ is a pre-seed under~$\mu$, then
$\theta(x)$ is a pre-seed under $\theta(\mu)$. Therefore if
$C$ is an equivalence class of pre-seeds under $\mu\in\XX$,
then $\theta(C)$ is an equivalence class of pre-seeds under~$\theta(\mu)$.
Also, if $c \in\R^d$ is the center of the ball with
the smallest radius that contains $C$, then $\theta(c)$ is the center
of the ball with the smallest radius that contains $\theta(C)$. Hence
if $b$
is a globe under $\mu$, then $\theta(b)$ is a globe under~$\theta(\mu)$. So clearly, $\Psi$ is isometry-equivariant.

Property \hyperlink{good}{(c)}: Let $\ww, \ww'\,{\in}\,\XX$, and
assume that $\ww\,{=}\,\ww'$ on $H(\ww)$. By Lemma~\ref{samepreseed}, $\ww$ and $\ww'$ have the same pre-seeds.
Thus, they have the same seeds, and hence the same globes.
Therefore by the definition of $\Psi$, we have $\Psi(\ww) = \Psi(\ww')$.
\end{pf*}

\section{Encoding and distributing randomness}
\label{pretheorem} Unfortunately, our proofs of Theorems
\ref{result} and \ref{resulttwo} do not follow from
Proposition \ref{fund} alone. Recall that in Examples
\ref{firstfirsteg} and \ref{firsteg} we partitioned $\R^d$
into cubes, and the cubes that contained exactly one or two
Poisson points were special. The locations of the Poisson
points in a special cube were converted into sequences of
i.i.d. $\uniforma$ random variables whose elements were then
assigned to the other cubes of the partition. The purpose of
this section is to state a lemma that asserts the existence of
a function that encapsulates the task of encoding and
distributing randomness in the more complicated case where a
deterministic partition is replaced by the selection rule from
Section \ref{goodclumping}, and Example \ref{firsteg} is
replaced by Theorem \ref{result}.

Let $\Psi$ be a selection rule. We say that a globe under
$\mu$ is \textit{one-special} if it happens to contain exactly
one $\mu$-point, and \textit{two-special} if it happens to
contain exactly two $\mu$-points. A globe is \textit{special} if
it is either one-special or two-special. Denote the set of
one-special globes by $\oneglobe[\Psi(\mu)]$, the set of two-special
globes by $\twoglobe[\Psi(\mu)]$ and the set of special globes by
$\sglobe[\Psi(\mu)]$. Also let
$\onesetglobe(\mu),\twosetglobe(\mu)$ and $\ssetglobe(\mu)$
denote the union of the set of one-special, two-special and
special globes, respectively. Let $\outsideonemu,
\outsidetwomu$ and $\outsidemu$ denote the respective
complements in $\R^d$. Note that by Proposition \ref{fund},
if $X$ is a Poisson point process on $\R^d$ with positive
intensity and law $\law$, then one-special globes and
two-special globes exist under $X$ $\law$-a.s.
\begin{lemma}[(Assignment function)]
\label{existenceassign} Let $d \geq1$ and $R > 0$. Let
$\Psi= \Psi_{R}$ be the selection rule from Section
\ref{goodclumping}. There exists a function
$\ass=\ass_{\Psi}: \XX\times(\borelf\cup\R^d) \to[0,1]$
with the following properties.
\begin{longlist}[(a)]
\item[(a)]
\hypertarget{assiid} Let $X$ be a Poisson point process on
$\R^d$ with positive intensity. Let $\{ \order(X)_i \}_{i
\in\N} := \operatorname{Globes}[\Psi(X)] \cup[X\ronnn_{\Psi(X)
^c} ]$,
where we have ordered the set using the radial ordering.
(Recall that globes are ordered by their centers.) If
$\{ U_i \}_{i \in\N}$ is a sequence of i.i.d. $\uniforma$
random variables that is independent of $X$, then
%
%
\begin{eqnarray}
\label{assiidstatement}
&&
\bigl( X \ronnn_{\outsideone}, \onesetglobe(X), \Psi(X),
\{ \ass(X, \order(X)_i) \}_{i \in\N} \bigr) \nonumber\\[-8pt]\\[-8pt]
&&\qquad\eqd \bigl( X \ronnn_{\outsideone}, \onesetglobe(X), \Psi(X),
\{ U_i \}_{i \in\N} \bigr)\nonumber
\end{eqnarray}
and
%
%
\begin{eqnarray}
\label{assiidstatementtwo}
&&\bigl( X \ronnn_{\outsidetwo}, \twosetglobe(X), \Psi(X),
\{ \ass(X, \order(X)_i) \}_{i \in\N} \bigr) \nonumber\\[-8pt]\\[-8pt]
&&\qquad\eqd \bigl( X \ronnn_{\outsidetwo}, \twosetglobe(X),
\Psi(X), \{ U_i \}_{i \in\N} \bigr) .\nonumber
\end{eqnarray}
\item[(b)]
\hypertarget{assiso} The map $\ass$ is isometry-invariant;
that is, for all isometries $\theta$ of $\R^d$ and for
all $(\mu, b) \in\XX\times(\borelf\cup\R^d)$, we
have $\ass(\mu, b) = \ass(\theta(\mu), \theta(b))$.
\end{longlist}
\end{lemma}

We call $\mathbf{U}_{\Psi}$ the \textit{assignment function}
for the selection rule $\Psi$. Thus if $X$ is a Poisson point
process and $b \in\operatorname{Globes}[\Psi(X)]$ or if $b \in
[X\ronnn_{\Psi(X) ^c}]$, then the assignment function
assigns a $\uniforma$ random variable $\ass(X, b)$ to $b$.
Property~\hyperlink{assiid}{(a)} states that the $\uniforma$ random
variables have a certain independence property; the values of
$X$ on both the one-special and two-special globes are needed
to determine the values of the assignment function. The map
that we shall define in the next section to prove Theorem
\ref{result} will use $\ass$ to assign $\uniforma$ random
variables to the globes and the points of the ether. We shall
see that property \hyperlink{assiid}{(a)} makes proving Theorem
\ref{result} easy. Property \hyperlink{assiso}{(b)} is necessary to
ensure that the map that we define is isometry-equivariant.

Let us also remark that since by property \hyperlink{good}{(c)} in the
definition of a selection rule, $\Psi(X)$ depends only on
$X\ronnn_{\Psi(X)^c} \subset X\ronnn_{\outside}$; therefore the
addition of $\Psi(X)$ in (\ref{assiidstatement}) and
(\ref{assiidstatementtwo}) is actually redundant. We now
have all the tools we need to prove Theorem \ref{result}.
We defer the proof of Lemma \ref{existenceassign} to Section~\ref{deferproof}.
Much of the proof is bookkeeping, but for
property \hyperlink{assiid}{(a)} we shall need to appeal to Proposition
\ref{fund}.

\section{\texorpdfstring{Proof of Theorem \protect\ref{result}}{Proof of Theorem 1}}
\label{def}
We are now in a position to prove Theorem~\ref{result}.
First we give the definition of the mapping that satisfies the
conditions of Theorem \ref{result}. Let $X$ be a Poisson
point process on $\R^d$ with intensity $\lx$, and let $\ly<
\lx$. Recall the definition of the splitting $\gammathin$ from
Section \ref{examples} (Proposition~\ref{specquant}) and the
definitions of $\relstandthin$ and $\standthin$ from
(\ref{coin}) and (\ref{standthin}) of Section
\ref{goodclumping}. Let $R= R(\lx, \ly) > 0$ be so that the
Lebesgue measure of $B(\bfz,R)$ is larger than the constant
$K(\lx, \ly)$ of Proposition~\ref{quantile}. Let $\Psi=
\Psi_{R}$ be the $R$-selection rule from Section
\ref{goodclumping}, and let $\ass$ be the assignment function
from Lemma \ref{existenceassign}. Define $\Phii=
\Phii_{(\lx, \ly)}$ as follows. For all $\mu\in\XX$,
%
%
\begin{eqnarray}
\label{Phii}
\Phii(\mu)&:=& \sum_{ b \in\operatorname{Globes}[\Psi(\mu)]}
\gammathin_{(b, \lx, \ly)}(\mu\ronnn_b, \ass(\mu, b))\nonumber\\[-8pt]\\[-8pt]
&&{} +
\sum_{x \in[ \mu\ronnn_{\Psi(\mu) ^c}]}
\relstandthin_{(\lx, \ly)}(x,\ass(\mu, {x})) .\nonumber
\end{eqnarray}
\begin{pf*}{Proof of Theorem \ref{result}}
From the definition of $\Phii$ it is easy to check that it is
isometry-equivariant; we need only recall that by Lemma
\ref{existenceassign}, the assignment function $\ass$ is
isometry-invariant and that the splitting $\gammathin$ and
selection rule $\Psi$ are isometry-equivariant. Also it is
obvious that $\Phii$ is monotone, so it suffices to check that
$\Phii(X)$ and $X - \Phii(X)$ are Poisson point processes on
$\R^d$ with intensities $\ly$ and $\lx- \ly$, respectively.

Let $U$ be a $\uniforma$ random variable independent of $X$, and
let $g_i\dvtx[0,1] \to[0,1]$ be the functions from Lemma
\ref{reprod}. Let $\{ b_i \}_{i \in\Z^{+}} = \operatorname
{Globes}[\Psi(X)]$,
where we have ordered the globes via the radial ordering.
Similarly, let $\{ x_i \}_{i \in\Z^{+}}= [X \ronnn_{\Psi(X)
^c}]$. Let $\Phi$ be the splitting defined
in Corollary \ref{fresh}. Note that $\Phi$ is a~version of
$\Phii$ that uses randomness from $U$ instead of from certain
points of~$X$. By property \hyperlink{interest}{(d)} of Proposition
\ref{specquant} we have that
%
%
\begin{eqnarray}
\label{phiversiona}
\Phi(X, U) &=& \sum_{i \in\Z^{+}} \ind_{[X(b_i) \not= 1]}
\gammathin_{b_i}(X \ronnn_{b_i}, g_i(U)) \nonumber\\[-8pt]\\[-8pt]
&&{} + \standthin\bigl(X \ronnn_{\Psi(X) ^c},
g_0(U)\bigr) \nonumber
\end{eqnarray}
and
\begin{eqnarray}
\label{phiversionb}
X - \Phi(X, U)
&=& \sum_{i \in\Z^{+}} \ind_{[X(b_i) \not= 2]}
\bigl( X\ronnn_{b_i} - \gammathin_{b_i}(X \ronnn_{b_i},
g_i(U)) \bigr)
\nonumber\\[-8pt]\\[-8pt]
&&{} +
X \ronnn_{\Psi(X) ^c} - \standthin\bigl(X \ronnn_{\Psi(X)
^c}, g_0(U)\bigr) .\nonumber
\end{eqnarray}
We shall show that $\Phii(X) \eqd\Phi(X,U)$ and $X - \Phii(X) \eqd
X - \Phi(X, U)$.
Set
\begin{eqnarray*}
\alpha&:=& \sum_{i \in\Z^{+}} \gammathin_{b_i}(X \ronnn_{b_i},
\ass(X, b_i)) , \\
\beta&:=& \sum_{i \in\Z^{+}} \relstandthin(x_i, \ass(X, x_i)) ,
\\
\alpha' &:=& \sum_{i \in\Z^{+}} \bigl(X\ronnn_{b_i} - \gammathin_{b_i}
(X \ronnn_{b_i}, \ass(X, b_i))\bigr)
\end{eqnarray*}
and
\[
\beta' :=
X\ronnn_{\Psi(X) ^c} - \sum_{i \in\Z^{+}} \relstandthin(x_i, \ass
(X, x_i)) .
\]
By definition,
\[
\Phii(X) = \alpha+ \beta
\quad\mbox{and}\quad X - \Phii(X) = \alpha' + \beta' .
\]
By property \hyperlink{interest}{(d)} of Proposition \ref{specquant},
\[
\alpha = \sum_{i \in\Z^{+}}\ind_{[X(b_i) \not= 1]}
\gammathin_{b_i} (X \ronnn_{b_i}, \ass(X, b_i))
\]
and
\[
\alpha' = \sum_{i \in\Z^{+}}\ind_{[X(b_i) \not= 2]}
\bigl(X \ronnn_{b_i} - \gammathin_{b_i} (X \ronnn_{b_i}, \ass(X,
b_i))\bigr) .
\]
By property \hyperlink{assiid}{(a)} of Lemma \ref{existenceassign}, we have that
\[
\alpha+ \beta \eqd \sum_{i \in\Z^{+}} \ind_{[X(b_i) \not=1 ]}
\gammathin_{b_i}(X \ronnn_{b_i}, g_{2i}(U)) + \sum_{i \in\Z^{+}}
\relstandthin(x_i, g_{2i + 1}(U))
\]
and
\begin{eqnarray*}
\alpha' + \beta' &\eqd& \sum_{i \in\Z^{+}}
\ind_{[X(b_i) \not=2 ]} \bigl( X \ronnn_{b_i} - \gammathin_{b_i}
(X \ronnn_{b_i}, g_{2i}(U)) \bigr)\\
&&{} + X\ronnn_{\Psi(X) ^ c} - \sum_{i \in\Z^{+}}
\relstandthin(x_i, g_{2i + 1}(U)) .
\end{eqnarray*}

Thus, by (\ref{phiversiona}), (\ref{phiversionb}), (\ref{coin}) and
(\ref{standthin}) it is easy to see that
\[
\alpha+ \beta \eqd \Phi(X, U)
\]
and
\[
\alpha' + \beta' \eqd X - \Phi(X, U) .
\]

Hence by Corollary \ref{fresh}, we have that $\Phii(X)$ and
$X - \Phii(X)$ are Poisson point processes with intensities
$\ly$ and $\lx- \ly$, respectively.
\end{pf*}

\section{The assignment function}
\label{deferproof}

In this section we shall prove Lemma \ref{existenceassign}.
Many of the same tools will be useful again in the proof of
Theorem \ref{resulttwo}. Recall that the assignment function\vadjust{\eject}
contained within it the two tasks of generating and
distributing uniform random variables. First we discuss how we
generate uniform random variables.

The following lemma describes explicitly how we convert the
position of a single $X$-point in a ball (which is a uniform
random variable on the ball) into a single uniform random
variable on $[0,1]$. We need to be explicit to preserve
equivariance.
\begin{lemma}[(Uniform random variables)]
\label{uniform}
For every\vspace*{1pt} $d \geq1$, $c \in\R^d$ and $R > 0$, define $f_{\bar
B(c,R)}\dvtx
\bar{B}(c, R) \to[0,1]$ via
\[
f_{\bar{B}(c,R)}(x) :=
\biggl(\frac{\|x-c\|}{R}\biggr)^d .
\]
The collection of
mappings $\{ f_{\bar{B}(c, R)} \}_{c \in\R^d}$ has the
following properties:
\begin{enumerate}
\item If $V$ is a $\uniformb[\bar{B}(c, R)]$ random variable,
then $f_{\bar{B}(c,R)}(V)$ is a $\uniforma$ random variable.
\item
\hypertarget{isof} We have isometry-invariance; that is,
for
any isometry $\theta$ of $\R^d$,\break
$f_{\bar{B}(c,R)}(x)=f_{\theta(\bar{B}(c,R))}(\theta(x))$ for
all $x \in\bar{B}(c,R)$.
\end{enumerate}
\end{lemma}
\begin{pf}
Here we shall make good use of the fact that we are working
with balls. Recall that the Lebesgue measure of a $d$-ball of
radius $R$ is given by $C(d) R^d$ for some fixed constant $C(d)
> 0$ depending only on $d$. Let $V$ be a uniform random
variable on the ball $\bar{B}(c,R)$. Then for $0 \le x \le1$,
\[
\PP\bigl(f_{\bar{B}(c,R)}(V) \leq x\bigr) = \PP(\|V - c\| \leq Rx^{{1/d}})
= \frac{ \les(\bar{B}(Rx^{{1/d}} )) }{ \les(\bar{B}(R)) } = x.
\]
\upqed\end{pf}

Each globe or $X$-point not in a globe will be associated to a
one-special globe and to a two-special globe. It will be
necessary to allow more than one globe or $X$-point to be
associated to each special globe. First we need to develop some
infrastructure. Recall that by Lemma \ref{reprod} a single
uniform random variable can be used to generate a sequence of
i.i.d. $\uniforma$ random variables.

\subsection*{Encoding functions}
We associate to every special globe a $[0,1]$-valued sequence
in the following way. Let $\{ f_{\bar{B}(c,R)} \}_{c \in\R^d}$
and $\{ g_i \}_{ i \in\N}$ be the collections of functions from
Lemmas \ref{uniform} and \ref{reprod}, respectively. Let
$\Psi$ be an $R$-selection rule. For each $b \in
\oneglobe[\Psi(\mu)]$, let $x_b$ denote the unique $\mu$-point in
$b$, and for each $b \in\twoglobe[\Psi(\mu)]$, let $x_b ^1$ and $x_b
^2$ be the two $\mu$-points in $b$, where we take $x_b ^1$ to
be the one closest to the origin in a lexicographic ordering.
Recall that $\modplus$ denotes addition modulo one. Let $\h=
\h_{\Psi}\dvtx\XX\times\goth{F} \to[0,1]$ and $h= h_{\Psi}\dvtx
\XX\times\goth{F} \to[0,1]^{\N}$ be defined as follows:
%
%
\begin{eqnarray}
\label{simpencodingfunction}\quad
\h_{\Psi}(\mu, b) &:=& \cases{
f_b(x_b), &\quad if $\mu\in\XX, b \in\oneglobe[\Psi(\mu)]$,\cr
f_b(x_b^1) \modplus f_b(x_b^2), &\quad if $\mu\in\XX, b \in
\twoglobe[\Psi(\mu)]$,\cr
0, &\quad if $\mu\in\XX, b \notin\sglobe[\Psi(\mu)]$,}
\end{eqnarray}
and
%
%
\begin{equation}
\label{encodingfunction}
h_{\Psi}(\mu, b) := \{ g_i(\h(\mu, b)) \}_{i \in\N} .
\end{equation}
We call $h_{\Psi}$ the \textit{encoding function} for the
selection rule $\Psi$, and we call $\h_{\Psi}$ the
\textit{simplified encoding function} for the selection rule
$\Psi$.
\begin{lemma}
\label{sequences} Let $d \geq1$ and $R > 0$. Let $\Psi$ be
an $R$-selection rule. Both the encoding and the simplified
encoding functions $h, \h$ satisfy the following properties:
\begin{longlist}[(a)]
\item[(a)] The maps $h,\h$ are isometry-invariant; that is,
for all isometries $\theta$ of~$\R^d$ and for all
$(\mu, b) \in\XX\times\goth{F}$, $h(\mu, b) =
h(\theta(\mu), \theta(b))$ and $\h(\mu, b) =
\h(\theta(\mu),\allowbreak \theta(b))$.
\item[(b)]
\hypertarget{jointdist} Let $X$ be a Poisson point process on
$\R^d$ with positive intensity. Let $\{ b^1_i \}_{i \in\N}:=
\oneglobe[\Psi(X)]$ and $\{ b^2_i \}_{i \in\N}:= \twoglobe[\Psi(X)]$,
where we have ordered the sets of one-special and two-special
globes by the radial ordering. If $\{ U_i \}_{i \in\N}$ is
a sequence of i.i.d. $\uniforma$ random variables that is
independent of $X$, then
\[
\bigl(X \ronnn_{\outsideone}, \onesetglobe(X), \{ \h(X,b^1_i) \}
_{i \in\N} \bigr)
\eqd\bigl(X \ronnn_{\outsideone},\onesetglobe(X), \{ U_i \}_{i
\in
\N} \bigr)
\]
and
\[
\bigl(X \ronnn_{\outsidetwo}, \twosetglobe(X), \{ \h(X,b^2_i) \}
_{i \in\N} \bigr)
\eqd\bigl(X \ronnn_{\outsidetwo},\twosetglobe(X), \{ U_i \}_{i
\in
\N} \bigr) .
\]
Similarly, if $\{ U'_i \}_{ i \in\N}$ is a sequence of i.i.d.
random variables, independent of $X$, where $U'_1 \eqd\{ U_i \}_{i
\in
\N}$, then
\[
\bigl(X \ronnn_{\outsideone}, \onesetglobe(X), \{ h(X,b^1_i) \}_{i
\in\N} \bigr)\eqd\bigl(X \ronnn_{\outsideone},
\onesetglobe(X), \{ U'_i \}_{i \in\N} \bigr)
\]
and
\[
\bigl(X
\ronnn_{\outsidetwo}, \twosetglobe(X), \{ h(X,b_i ^2) \}_{i \in\N}
\bigr) \eqd\bigl(X \ronnn_{\outsidetwo}, \twosetglobe(X),
\{ U'_i \}_{i \in\N} \bigr) .
\]
\end{longlist}
\end{lemma}
\begin{pf}
The proof of property (a) follows immediately from the
definition of an encoding function, property \hyperlink{seliso}{(b)} of
a selection rule and Lem\-ma~\ref{uniform}. We now focus our
attention on property \hyperlink{jointdist}{(b)}. From the definition
of $h$ and the fact that the $g_i$ satisfy the conditions of
Lemma \ref{reprod}, it suffices to verify the condition for
the simplified encoding function $\h$.

We need some additional notation. Let $\{ b_i \}_{i \in\N} =
\operatorname{Globes}[\Psi(X)]$,
where we have ordered the globes via the radial ordering. Let $c_i \in
b_i$ be the centers of the globes. Also let $c^1 _i \in b ^1 _i$ and
$c^2 _i \in b^2 _i$ be the centers of the one-special and two-special
globes. Let $\gothh{V^1_i}$ be the unique $X$-point in each $b^1_i$.
Similarly, let~$\gothh{V^2_i}$ be the set of unordered
$X$-points of each $b^2_i$. Let $\theta_{y}$ be the isometry of~$\R
^d$ such that for all
$x \in\R^d$, we have $\theta_y(x) = x +y$. Assume $X$ has intensity~$\lx$.

It follows from Proposition \ref{fund} that $ \{ \theta^{-1}_{c_i}
(X \ronnn_{b_i}) \}_{i \in\N}$ is a sequence of i.i.d. Poisson point
processes on $\bar{B}(\bfz, R)$ with intensity $\lx$; furthermore,
the sequence is independent of $(X \ronnn_{\Psi(X) ^c}, \Psi
(X))$. Hence by Lemma \ref{char}, $\{ \gothh{V^1_i} - c^1_i  \}
_{i \in\N}$ is a sequence of i.i.d. $\uniformb[\bar{B}(\bfz,R)]$
random variables that is
independent of $(X \ronnn_{ \outsideone}, \onesetglobe(X))$.
Similarly, we
have that $\{ \theta^{-1}_{c^2 _i}\gothh{V^2_i}  \}_{i \in\N}$ is a
sequence\vspace*{-2pt} of i.i.d. pairs of unordered
$\uniformb[\bar{B}(\bfz,R)]$ random variables that is
independent of $( X\ronnn_{\outsidetwo}, \twosetglobe(X))$.
By the definition of $\h$ and Lemmas \ref{uniform} and
\ref{modulo}, the result follows immediately.
\end{pf}

We turn now to the task of distributing randomness. A natural
approach is to have each nonspecial globe request randomness
from the closest available special globe (where distances are
measured between the centers of the globes). However, we do
not know much about the process of globe-centers. In
particular, it is not immediately obvious that it has distinct
inter-point distances $\law$-a.s. To avoid this
problem, we shall make use of some of the other properties of
seeds. Recall that if $x$ is a pre-seed, we call $A(x;R+10, 3R +75 +
d)$ the halo. If $x$ is a seed, then we shall also call
$A(x;R+10, 3R + 75 + d)$ the \textit{halo}. We shall associate to
every globe a point in its halo in an equivariant way.

\subsection*{Tags}
Let $\Psi$ be the selection rule from Section \ref{goodclumping},
and let the inter-point distances of $\mu\in\XX$ be distinct.
For each globe under $\mu$, we choose a point in its halo in
the following isometry-equivariant way. First note that the
halo contains more than three $\mu$-points. Take the two
mutually closest points in the halo, then choose the one of
this pair that is closest to the other points in the halo. We
call this point the \textit{tag} of the globe. We note that by
Lemma \ref{interinter}, part \hyperlink{inter}{(b)}, tags are well
defined and exist for every globe $\law$-a.s. For
completeness, if the inter-point distances in the halo are not
distinct, we take the tag to be the center of the globe. Let
$t =t_{\Psi}\dvtx\XX\times\goth{F} \to\R^d \cup\{ \infty \}$ be
the measurable function defined as follows:
%
%
\begin{eqnarray}
\label{associatedpoints}
t_{\Psi}(\mu, b) &:=& \cases{
\mbox{the tag of }b, &\quad if $\mu\in\XX, b \in\operatorname
{Globes}[\Psi(\mu)]$,\cr
x, &\quad if $\mu\in\XX, b=\{ x \}, x \in\bigl[\mu\ronnn_{\Psi(\mu)
^c}\bigr]$,\cr
\infty, &\quad otherwise.}
\end{eqnarray}
We call $t_{\Psi}$ the \textit{tagging function} for the selection rule
$\Psi$.
\begin{lemma}
\label{tag} Let $d \geq1$ and $R > 0$. Let $\Psi=\Psi_{R}$
be the selection rule from Section \ref{goodclumping}. The
tagging function $t=t_{\Psi}\dvtx\XX\times\goth{F} \to\R^d \cup
\{ \infty \}$ has the following properties:
\begin{enumerate}
\item The map $t$ depends only on $(\Psi(\mu), \mu
\ronn_{\Psi(\mu) ^c})$; that is, for all $\mu, \mu' \in\XX$ if
$(\Psi(\mu), \mu\ronn_{\Psi(\mu) ^c}) =(\Psi(\mu'), \mu'
\ronn_{\Psi(\mu') ^c})$, then $t(\mu, \cdot)= t(\mu', \cdot)$.
\item The map $t$ is isometry-equivariant; that is, for
all isometries $\theta$ of $\R^d$ and for all $(\mu,
b) \in\XX\times\goth{F}$, $\theta(t(\mu, b)) =
t(\theta(\mu), \theta(b))$. Here we take
$\theta(\infty) = \infty$.
\end{enumerate}
\end{lemma}
\begin{pf}
The result follows immediately from the definition of the tagging function.
\end{pf}

\subsection*{Partners and ranks}
Let $\Psi$ be the selection rule from Section \ref{goodclumping}.
We shall now measure distances between globes, as well as
distances between globes and $\mu$-points, via the distances
between their tags. Let the inter-point distances of $\mu\in
\XX$ be distinct and also assume that $\oneglobe[\Psi(\mu)]$ and
$\twoglobe[\Psi(\mu)]$ are both nonempty. For each globe $b \in
\operatorname{Globes}[\Psi(\mu)]$ we call its closest one-special
globe its
\textit{one-partner}, and its closest two-special
globe its \textit{two-partner}. Similarly, for each
$x \in[\mu\ronnn_{\Psi(\mu) ^c }]$ we call its closest
one-special globe its \textit{one-partner} and its closest
two-special globe its \textit{two-partner}. Suppose that a
globe $b$ has a special globe $B \in\sglobe[\Psi(\mu)]$ as a
partner; then $B$ assigns the number $2n$ to $b$ if there are
exactly $n$ globes with $B$ as partner that are closer to $B$
than $b$. We call the number that $b$ is assigned by its
one-partner its \textit{one-rank} and the number that $b$ is
assigned by its two-partner its \textit{two-rank}. Similarly,
a special globe $B \in\sglobe[\Psi(\mu)]$ assigns the number $2n+1$
to $x$ if it is a partner of $x \in[\mu\ronnn_{\Psi(\mu) ^c}]$,
and there are exactly~$n$ partners in $[\mu]$ that
are closer to $B$ than $x$; we also call the number that $x$ is
assigned its \textit{one-rank} or \textit{two-rank} depending on
whether it is assigned by its one- or two-partner. Let $\XXd$
be the set of point measures of~$\XX$ that have both one- and
two-special globes and have distinct inter-point distances.
We define $p = p_{\Psi}\dvtx\XX\times(\goth{F} \cup\R^d) \to
\goth{F} \times\goth{F}$ as follows:
\[
p_{\Psi}(\mu, b) :=
(\mbox{one-partner of }b, \mbox{two-partner of }b)
\]
if $\mu\in\XXd$ and $b \in\operatorname{Globes}[\Psi(\mu)] \cup
[\mu\ronnn_{\Psi(\mu) ^c} ]$
and $p_{\Psi}(\mu, b) := (b,b)$ otherwise.
We also define $r = r_{\Psi}\dvtx\XX\times(\goth{F} \cup\R^d)
\to\N\times\N$ as follows:
\[
r_{\Psi}(\mu, b) :=
(\mbox{one-rank of }b, \mbox{two-rank of }b)
\]
if $\mu\in\XXd$ and $b \in\operatorname{Globes}[\Psi(\mu)] \cup
[\mu\ronnn_{\Psi(\mu) ^c} ]$
and $r_{\Psi}(\mu, b) := (0,0)$ otherwise.
We call $p_{\Psi}$ the \textit{partner function} for the
selection rule $\Psi$, and we call $r_{\Psi}$ the \textit{rank
function} for $\Psi$. Also let
\[
\chi(\mu) := \bigl(\Psi(\mu),
\onesetglobe(\mu), \twosetglobe(\mu), \mu\ronnn_{\Psi(\mu)
^c}\bigr)
\]
for all $\mu\in\XX$.
\begin{lemma}
\label{ordering} Let $d \geq1$ and $R > 0$. Let $\Psi$ be
the selection rule from Section~\ref{goodclumping}. The
partner and rank functions $p=p_{\Psi}$ and $r= r_{\Psi}$ have
the following properties:
\begin{enumerate}
\item The maps $p,r$ depend only on $\chi(\mu)$; that is, for
all $\mu, \mu' \in\XX$ if $\chi(\mu) = \chi(\mu')$, then
$p(\mu, \cdot)= p(\mu', \cdot)$ and $r(\mu, \cdot)= r(\mu',
\cdot)$.
\item The map $p$ is isometry-equivariant; that is, for
all isometries $\theta$ of $\R^d$ and for all $(\mu,
b) \in\XX\times\goth{F}$, $\theta(p(\mu, b)) =
p(\theta(\mu), \theta(b))$.
\item The map $r$ is isometry-invariant; that is, for
all isometries $\theta$ of $\R^d$ and for all $(\mu,
b) \in\XX\times\goth{F}$, $r(\mu, b) =
r(\theta(\mu), \theta(b))$.
\end{enumerate}
\end{lemma}
\begin{pf}
The result follows immediately from the definitions of the
partner and rank functions and Lemma \ref{tag}.
\end{pf}
\subsection*{Assignment functions}
We shall now combine the encoding, partner and rank functions
to obtain an assignment function. Let $\Psi$ be the selection
rule from Section \ref{goodclumping}. Define $\mathbf{U}
= \mathbf{U}_{\Psi}\dvtx\XX\times(\borelf\cup\R^d) \to
[0,1]$ as follows. Let $h=h_{\Psi}$, $p=p_{\Psi}$ and
$r=r_{\Psi}$ be the encoding, partner and rank functions.
Recall that $h\dvtx\XX\times\borelf\to[0,1] ^{\N}$. For all
$(\mu,b) \in\XX\times(\borelf\cup\R^d)$, let
%
%
\begin{equation}
\label{assignmentfunction}
\mathbf{U}(\mu, b) := h(\mu, p(\mu, b)_1)_{r(\mu, b)_1}
\modplus h(\mu, p(\mu, b)_2)_{r(\mu, b)_2} .
\end{equation}
\begin{pf*}{Proof of Lemma \ref{existenceassign}}
The isometry-invariance of $\ass$ follows immediately
from the definition of $\ass$ and Lemmas \ref{sequences} and
\ref{ordering}. Let $X$ be a Poisson point process in
$\R^d$. Let $\{ \order_i \}_{i \in\N}:= \operatorname
{Globes}[\Psi(X)] \cup
[X \ronnn_{\Psi(X) ^c}]$, $\{ b^1_i \}_{i \in\N}:=
\oneglobe[\Psi(X)]$, and let $\{ b^2_i \}_{i \in\N}:=
\twoglobe[\Psi(X)]$, where we have ordered the sets using the radial
ordering. Let $\{ U_i \}_{ i \in\N}$ be a sequence of i.i.d.
$\uniforma$ random variables independent of $X$. Let
$\{ U_i' \}_{i \in\N}$ be an i.i.d. sequence (independent of
$X$), where $U_1'$ is a sequence of i.i.d. $\uniforma$ random
variables. From Lemma \ref{ordering}, $p(X, \cdot)$ and
$r(X, \cdot)$ depend only on $\chi(X)$. It is clear that both
$\chi(X)$ and $h(X, b^2_i)$ depend only on
$(X\ronnn_{\outsideone}, \onesetglobe(X))$, so that by Lemma
\ref{sequences}
%
%
\begin{eqnarray}
\label{toexists}
&&\bigl(X \ronnn_{\outsideone}, \onesetglobe(X), \chi(X), \{
h(X,b^1_i) \}_{i \in\N},
\{ h(X,b^2_i) \}_{i \in\N} \bigr) \nonumber\\[-8pt]\\[-8pt]
&&\qquad\eqd \bigl(X \ronnn_{\outsideone}, \onesetglobe(X), \chi(X), \{ U'_i
\}_{i \in\N},
\{ h(X,b^2_i) \}_{i \in\N} \bigr) .\nonumber
\end{eqnarray}
From the definition of the assignment function, it is clear
that $\ass(X, \cdot)$ depends only on
\[
(\chi(X), \{ h(X, b^1_i) \}_{ i \in\N}, \{ h(X, b^2_i) \}
_{i \in\N}) .
\]
It
is also easy to see that $\chi(X)$ depends only on $(X
\ronnn_{\outside}, \onesetglobe(X), \twosetglobe(X))$. Thus
from the definition of the assignment function,
(\ref{toexists}) and Lemma \ref{modulo}, it follows that
\[
\bigl(X \ronnn_{\outsideone}, \onesetglobe(X), \{ \ass(X,\order
_i)\}_{i \in\N}
\bigr) \eqd\bigl(X \ronnn_{\outsideone}, \onesetglobe(X), \{ U_i
\}_{i \in\N} \bigr).
\]

Similarly, we have that
%
%
\begin{eqnarray}
\label{toexiststwo}
&&\bigl(X \ronnn_{\outsidetwo},  \twosetglobe(X), \chi(X), \{
h(X,b^1_i) \}_{i \in\N},
\{ h(X,b^2_i) \}_{i \in\N} \bigr) \nonumber\\[-8pt]\\[-8pt]
&&\qquad \eqd\bigl(X \ronnn_{\outsidetwo}, \twosetglobe(X), \chi(X),\{
h(X,b^1_i) \}_{i \in\N},
\{ U'_i \}_{i \in\N} \bigr) ,\nonumber
\end{eqnarray}
from which it follows that
\[
\bigl(X \ronnn_{\outsidetwo}, \twosetglobe(X), \{ \ass(X,\order
_i) \}_{i \in\N} \bigr)
\eqd\bigl(X \ronnn_{\outsidetwo}, \twosetglobe(X), \{ U_i \}_{i
\in
\N} \bigr) .
\]
\upqed\end{pf*}

%
\section{\texorpdfstring{Proof of Theorem \protect\ref{resulttwo}}{Proof of Theorem 2}}
\label{proofoftheoremtwo} In this section, we shall show how
the tools used to prove Theorem \ref{result} can be adapted
to prove Theorem \ref{resulttwo}. As a first step we prove a source-universal
translation-equivariant version of Theorem~\ref{resulttwo}.
That is, given $\lambda'$, we define a
translation-equivariant map $\phithmtwo\dvtx\XX\to\XX$ such that
if $X$ is a Poisson process on $\R^d$ of any positive intensity
$\lambda$,
then $\phithmtwo(X)$ is a Poisson process of intensity
$\lambda'$. By modifying the map $\phithmtwo$ we shall obtain
a~map $\thmtwo$ that is isometry-equivariant and satisfies
the conditions of Theorem~\ref{resulttwo}. We need some
preliminary definitions before we can give the definition of~$\phithmtwo$.

\subsection*{Voronoi cells}
The Voronoi tessellation of a simple point measure $\mu\in
\XX$ is a partition of $\R^d$ defined in the following way. The
\textit{Voronoi cell} of a point $x \in[\mu]$ is the set of
all points $y \in\R^d$ such that $\|x-y\| < \|z-y\|$ for all
$z \in[\mu] \setminus\{ x \}$. The \textit{unclaimed}
points are the points that do not belong to a~cell. We define
the \textit{Voronoi tessellation} $\vp(\mu)$ to be the set of
all Voronoi cells along with the set of unclaimed points. Note
that if $\mu$ is locally finite and not identically zero, then
the set of unclaimed points has zero Lebesgue measure. Note
that the Voronoi tessellation is clearly isometry-equivariant; that is,
for any isometry $\theta$ of $\R^d$ we
have $\vp(\theta\mu) = \theta\vp(\mu): = \{ \theta\upsilon\dvtx
\upsilon
\in
\vp(\mu) \}$.

For each $A \in\borel$ with positive finite Lebesgue measure,
let $c_A$ be its center of mass. Let $\Psi$ be the $R$-selection
rule, and define $\cglobe=\cglobe_{\Psi}\dvtx\XX\to\XX$ via
%
%
\begin{equation}
\label{centerglobe}
\cglobe(\mu):= \sum_{b \in\operatorname{Globes}[\Psi(\mu)] }
\delta_{c_b} .
\end{equation}
Note that $\cglobe$ is also isometry-equivariant. The map
$\phithmtwo$ will be defined by placing independent Poisson
point processes in each Voronoi cell of $\vp(\cglobe(\mu))$. Recall
that the globes do not intersect so that a Voronoi cell will always
contain the globe with the same center. Let $\theta_{y}$ be the
isometry of $\R^d$ such that for all
$x \in\R^d$, we have $\theta_y(x) = x +y$. We define
$\phithmtwo$ in the following way. Let $\Psi$ be the
$R$-selection rule from Section \ref{goodclumping}, where we may
choose $R=1$. Let $\poiseed$ be the collection of
mappings from Lemma \ref{poiseed}. Let $\ass$ be the
assignment function from Lemma \ref{existenceassign}. The
map $\phithmtwo= \phithmtwo_{\ly}$ is defined via
%
%
\begin{equation}
\label{transdef}
\phithmtwo(\mu) := \sum_{ \upsilon\in\vp(\cglobe(\mu))}
\sum_{ b \in\operatorname{Globes}[\Psi(\mu)]} \ind_{[b \subset
v]}\theta_{c_\upsilon}
\bigl(\poiseed_{ ( \theta^{-1}_{c_\upsilon}(\upsilon) , \ly) } (\ass(\mu,
b)) \bigr) .
\end{equation}
Note that $\phithmtwo$ depends only on the parameter $\ly$, since
$\ass$ depends only on~$R$, which we have set equal to $1$.
\begin{proposition}
\label{phithmtwolabel} The map $\phithmtwo$ has the following
properties:
\begin{longlist}[(a)]
\item[(a)] The map $\phithmtwo$ is translation-equivariant.
\item[(b)]
If $X$ is a Poisson process on $\R^d$ with
positive intensity, then $\phithmtwo(X)$ is a Poisson process on $\R
^d$ with intensity $\ly$.

\end{longlist}
\end{proposition}
\begin{pf}
Part (a) follows from the fact that the assignment function is
isometry-invariant and that the selection rule, Voronoi
tessellation and the map $\cglobe$ are all isometry-equivariant. From
the definition of $\phithmtwo$, one can
verify that it is translation-equivariant since any two
translations of $\R^d$ commute with each other. However, since
translations and reflections do not necessarily commute, we
have only translation-equivariance. Part (b) follows
from Lemma \ref{poiseed} and Lemma \ref{existenceassign}
once we note that the Voronoi tessellation and the centers of
the globes and Voronoi cells depend only on $\Psi(X)$.
\end{pf}

The following example elaborates on the difficulty of defining
an isometry-equivariant version of $\phithmtwo$.
\begin{example} \label{no}
Let $\ly> 0$. Let $\borelbd\subset\borel$ be the set of
Borel sets with positive finite Lebesgue measure. There does
not exist a family of measurable functions $\poiseedd$ such
that for each $A \in\borelbd$, $\poiseedd_A\dvtx[0,1] \to\XX$
has the following properties:
\begin{enumerate}
\item If $U$ is a $\uniforma$ variable, then $\poiseedd_A(U)$
is a Poisson point process on $A$ with intensity $\ly$.
\item The map $\poiseedd$ is isometry-equivariant; that
is, for all isometries $\theta$ of $\R^d$,
$\poiseedd_{\theta A}(U) = \theta\poiseedd_A(U)$.
\end{enumerate}
\end{example}
\begin{pf}
Toward a contradiction, let $\poiseedd$ satisfy the above
properties. For each $x \in\R^d$, let $x_i$ be the
$i$th coordinate. Consider $A := B(\bfz,1)$, the unit
ball centered at the origin, and let $A':= \{ x \in A\dvtx x_1 > 0\}$.
Let $\theta$ be the reflection of the first coordinate;
that is, if $y = (y_1,\ldots, y_d)$ for some $y_i \in\R$, then
$\theta(y) = (-y_1, y_2,\ldots, y_d)$. Let $U$ be a
$\uniforma$ random variable. The event $E:=
\{ \poiseedd_A(U)(A) = \poiseedd_A(U)(A') = 1\}$ occurs with
nonzero probability. However, \mbox{$\theta(A) = A$}, so that
whenever $E$ occurs, $\poiseedd_{\theta A}(U) \not= \theta
\poiseedd_{ A }(U)$.
\end{pf}

Note that in the proof of Example \ref{no}, the
counterexample used a set $A$ that is invariant under rotations
and reflections. One would guess that the Voronoi cells of a
random process such as the centers of the special globes should
lack such symmetries. However, rather than dealing with the
symmetries of the Voronoi cells, we proceed as follows.

Let $X$ be a Poisson process on $\R^d$ with positive intensity,
and let $\Psi$ be the selection rule from Section
\ref{goodclumping}. Let $b \in\operatorname{Globes}[\Psi(X)]$ and
for simplicity
assume that its center is at the origin. From the definition
of a globe, there will always be at least $d$ points in the
halo of a globe. We shall choose $d$ points from the halo and
use them to associate an isometry to the globe. By choosing from the
halo of $b$ in an equivariant way $d$
points $\{ x_1,\ldots, x_d \}$ that are linearly independent, we shall
define an isometry $\theta$ with the following properties:
\begin{enumerate}
\item We have $\theta(\bfz) = \bfz\in\R^d$.
\item For all $i,j$ such that $1 \leq i < j\le d$, we have
$\theta(x_i)_j = 0 \in\R$; that is, the $j$th
coordinate of $\theta(x_i) \in\R^d$ is zero for $j > i$.
\item For all $i$ such that $1 \leq i \leq d$, we have
$\theta(x_i)_i > 0$.
\end{enumerate}
Selecting $d$ points from the halo of a globe is an easy
extension of the idea of a tag of a globe. Also to prove that
such an isometry exists and is unique, we appeal to the tools
of linear algebra, in particular the QR factorization lemma.
\subsection*{Notations and conventions}
To use the tools of linear algebra, it will be convenient to
identify elements of $\R^d$ with column vectors; that is, $\R^d
= \R^{d \times1}$. Given an isometry $\theta$ of $\R^d$ and a
matrix $A \in\R^{d \times d}$, we let $\theta(A) \in\R^{d
\times d}$ be the matrix obtained by applying $\theta$ to each
of the columns of $A$. Let $\vec{1} \in\R^{1 \times d}$ denote the
row vector with all ones in its entries. Thus given $c \in\R^d$, $c
\vec{1}$ is the $d \times d$ matrix where each of its columns is equal
to $c$. We also denote the identity matrix by
$I \in\R^{d \times d}$.
\subsection*{$d$-tags}
Let $\Psi$ be the selection rule from Section \ref{goodclumping},
and let the inter-point distances of $\mu\in\XX$ be distinct.\vadjust{\goodbreak}
The \textit{$d$-tag} of a globe $b \in\operatorname
{Globes}[\Psi(\mu)]$
is a matrix $A \in\R^{d \times d}$ defined inductively as
follows. The first column of the matrix is the tag of $b$.
Given that the $(i-1)$th column is already defined,
the $i$th column is the $\mu$-point in the halo of $b$
that is closest to the $(i-1)$th column and is not
equal to any of the first $i-1$ columns. For completeness, if
the inter-point distances in the halo are not distinct, we take
the $d$-tag to be the matrix in $\R^{d \times d}$ where each
column vector is the center of the globe. Let
${\bar{t}}={\bar{t}}_{\Psi}\dvtx
\XX\times\borelf\to\R^{d \times d} \cup\{ \infty \}$ be the
measurable function defined as follows:
%
%
\begin{eqnarray}
\label{associatedpointsd}
{\bar{t}}_{\Psi}(\mu, b) &:=& \cases{
\mbox{the $d$-tag of } b, &\quad if $\mu\in\XX\mbox{ and } b
\in\operatorname{Globes}[\Psi(\mu)]$,
\cr
\infty, &\quad otherwise.}
\end{eqnarray}
We call ${\bar{t}}_{\Psi}$ the \textit{$d$-tagging
function} for the selection rule $\Psi$.\vspace*{3pt}
\begin{lemma}
\label{tagd} Let $d \geq1$ and $R > 0$. Let $\Psi= \Psi_{R}$
be the selection rule from Section \ref{goodclumping}. The
$d$-tagging function ${\bar{t}}={\bar{t}}_{\Psi}\dvtx\XX\times\borelf
\to
\R^{d \times d} \cup\{ \infty \}$ has the following
properties:
\begin{enumerate}
\item The map ${\bar{t}}$ depends only on $(\Psi(\mu), \mu
\ronn_{\Psi(\mu) ^c})$; that is, for all $\mu, \mu' \in
\XX$, if
$(\Psi(\mu), \mu\ronn_{\Psi(\mu) ^c}) = (\Psi(\mu'), \mu'
\ronn_{\Psi(\mu') ^c})$, then ${\bar{t}}(\mu, \cdot)= {\bar{t}}(\mu',
\cdot)$.
\item The map ${\bar{t}}$ is isometry-equivariant; that is,
for all isometries $\theta$ of $\R^d$ and for all
$(\mu, b) \in\XX\times\borelf$, we have
$\theta({\bar{t}}(\mu, b)) = {\bar{t}}(\theta(\mu), \theta(b))$. We
take $\theta(\infty) = \infty$.
\end{enumerate}
\end{lemma}
\begin{pf}
The result follows immediately from the definition of the $d$-tagging function.
\end{pf}\vfill{\eject}

We note that the $d$-tag of a globe is almost surely a
nonsingular matrix by Lemma \ref{interinter}\hyperlink{span}{(c)}.
The following lemma allows us to associate an isometry to each
globe and its $d$-tag. Recall that every isometry $\theta$ of
$\R^d$ that fixes the origin can be identified with a unique
orthogonal matrix $Q \in\R^{d \times d};$ that is, there is a
unique matrix $Q$ such that $QQ^{T} = Q^{T}Q = I \in\R^{d
\times d}$ and $Qx = \theta(x)$ for all $x \in\R^d = \R^{d
\times1}$. For background, see \cite{MR681482}, Chapter 1.
\begin{lemma}[(QR factorization)]
\label{linalg} For all $d \geq1$, if $A \in\R^{d \times d}$
is a square matrix, then there exists an orthogonal matrix $Q
\in\R^{d \times d}$ and an upper triangular matrix $\Delta\in
\R^{d \times d}$ such that $A=Q\Delta$. Furthermore, if $A$ is
nonsingular, then the factorization is unique if we require
the diagonal entries of $\Delta$ to be positive.
\end{lemma}

For a proof, see, for example, \cite{MR1084815}, Section 2.6.

\subsection*{Upper triangular matrices and fixing isometries}
Let $\Psi$ be a selection rule from Section
\ref{goodclumping}, and let $b \in\operatorname{Globes}[\Psi(\mu
)]$. The
\textit{upper triangular matrix for} $b$ is the matrix $\Delta
\in\R^{d \times d}$ defined as follows. Let $c_b$ be the
center of the globe $b$. Let $A' \in\R^{d \times d}$ be the
$d$-tag for the globe $b$. Let $A:= A' - c_b\vec{1} $. If $A$ is
singular, then we take $\Delta= \bfz\in\R^{d \times d}$.
Otherwise, by Lemma \ref{linalg}, there exists a unique
factorization such that $A = Q\Delta$, where $Q \in\R^{d
\times d}$ is an orthogonal matrix, and $\Delta\in\R^{d \times
d}$ is an upper triangular matrix such that all its diagonal
entries are positive. When $A$ is nonsingular, we say that
the unique isometry~$\sigma$ such that $\sigma(c_b) = \bfz\in
\R^d$ and $\sigma(A') = \Delta$ is the \textit{fixing isometry}
for the globe~$b$.

Let $\bolds{\Delta} = \bolds{\Delta} _{\Psi}\dvtx\XX
\times\borelf\to\R^{d \times d}$ be
the measurable function defined as follows:
\[
\bolds{\Delta} _{\Psi}(\mu, b) :=
\mbox{the upper triangular matrix for } b
\]
if $\mu\in\XX$ and $b \in\operatorname{Globes}[\Psi(\mu)]$,
while $\bolds{\Delta} _{\Psi}(\mu, b) := I \in\R^{d \times
d}$ otherwise.
Let $\bfz\in(\R^d)^{\R^d}$ be the function that sends every
element of $\R^d$ to $\bfz\in\R^d$. The \textit{fixing isometry
function} $\thetaa= \thetaa_{\Psi}\dvtx\XX\times\borelf\to
(\R^d) ^{\R^d}$ for the selection rule $\Psi$ is defined as
follows:
\[
\thetaa_{\Psi}(\mu, b) :=
\mbox{the fixing isometry for the globe $b$},
\]
if $\mu\in\XX$, $b \in\operatorname{Globes}[\Psi(\mu)]$, and the
$d$-tag of $b$
is nonsingular, while
$\thetaa_{\Psi}(\mu, b) := \bfz\in(\R^d) ^{\R^d}$ otherwise.
\begin{lemma}
Let $d \geq1$ and $R > 0$. Let $\Psi= \Psi_{R}$ be the
selection rule from Section \ref{goodclumping}. The map
$\bolds{\Delta} = \bolds{\Delta} _{\Psi}\dvtx\XX\times
\borelf\to\R^{d \times d}$ and the
fixing isometry function $\thetaa= \thetaa_{\Psi}\dvtx\XX\times
\borelf\to(\R^d) ^{\R^d}$ have the following properties:
\label{orient}
\begin{enumerate}
\item The maps $ \bolds{\Delta} $ and $\thetaa$ depend only on
$(\Psi(\mu),
\mu\ronn_{\Psi(\mu) ^c})$; that is, for all $\mu, \mu' \in
\XX$ if $(\Psi(\mu), \mu\ronn_{\Psi(\mu) ^c}) = (\Psi(\mu'),
\mu' \ronn_{\Psi(\mu') ^c})$, then $ \bolds{\Delta} (\mu,
\cdot)= \bolds{\Delta} (\mu',
\cdot)$ and\break $\thetaa(\mu, \cdot)= \thetaa(\mu', \cdot)$.
\item The map $ \bolds{\Delta} $ is isometry-invariant; that
is, for
all isometries $\theta$ of $\R^d$ and for all $(\mu,
b) \in\XX\times\borelf$, we have $ \bolds{\Delta} (\mu, b) =
\bolds{\Delta} (\theta(\mu), \theta(b))$.
\end{enumerate}
\end{lemma}
\begin{pf}
The first property follows immediately from the definitions of
the maps and Lemma \ref{tagd}. We prove the second property
in the following way. Let $\theta$ be an isometry of $\R^d$.
Let $A' \in\R^{d \times d}$ be a nonsingular matrix, let $a'
\in\R^d$ and set $A := A' -a'\vec{1}$. Let $A = Q\Delta$ be the
unique QR factorization of
$A$, where all the diagonal entries of $\Delta$ are positive.
From the definition of the upper triangular matrix for a globe,
it suffices to show that for some orthogonal matrix $Q
^{\prime\prime}$, we have
\[
\theta(A') - \theta(a')\vec{1} =
Q^{\prime\prime}\Delta.
\]
Note that there exists an
orthogonal matrix $Q'$ and $c \in\R^d$ such that for all $x
\in\R^d = \R^{d \times1}$, we have $\theta(x) = Q'x +
c$. Observe that
\begin{eqnarray*}
\theta(A') - \theta(a')\vec{1} & = & Q'A' + c\vec{1} -(Q'a' + c\vec
{1}) \\
& = & Q'(A' - a'\vec{1})
=  Q'A \\
& = & (Q'Q)\Delta.
\end{eqnarray*}
\upqed
\vspace*{3pt}\end{pf}

We are now ready to give the definition of the mapping that
satisfies the conditions of Theorem \ref{resulttwo}.
Set $R=1$, and let $\Psi$ be the
$R$-selection rule from Section \ref{goodclumping} with $R=1$.
Let $\thetaa\dvtx\XX\times\borelf\to(\R^d) ^{\R^d}$ be the
fixing isometry function for $\Psi$, let $\poiseed$ be a
collection of functions from Lemma \ref{poiseed} and let
$\ass$ be the assignment function from Lemma
\ref{existenceassign}.
Define $\thmtwo= \thmtwo_{\ly} \dvtx\XX\to\XX$ as
%
%
\begin{eqnarray}
\label{defthmtwo}
\thmtwo(\mu) &:=& \sum_{\upsilon\in\vp(\cglobe(\mu))}
\sum_{b \in\operatorname{Globes}[\Psi(\mu)]} \ind_{[b \subset v]}
\ind_{[\thetaa(\mu, b) \not= \bfz]}\nonumber\\[-8pt]\\[-8pt]
&&\hphantom{\sum_{\upsilon\in\vp(\cglobe(\mu))}
\sum_{b \in\operatorname{Globes}[\Psi(\mu)]} }
{}\times
\thetaa(\mu, b)^{-1} \bigl( \poiseed_{ (\thetaa(\mu, b)(\upsilon), \ly
) }
(\ass(\mu,b)) \bigr)\nonumber
\end{eqnarray}
for all $\mu\in\XX$.
\begin{pf*}{Proof of Theorem \ref{resulttwo}}
From the definition of $\thmtwo$ it is almost immediate that it
is isometry-equivariant. It suffices to check the following
claim. Let $\upsilon\in\borel$, let $b \in\operatorname{Globes}[\Psi
(\mu)]$, and let
$\theta$ be any isometry of $\R^d$. We claim that for all
$\mu\in\XX$,
%
%
\begin{eqnarray}
\label{claimiso}
&&\thetaa( \theta\mu, \theta b) ^{-1}
\bigl( \poiseed_{\thetaa( \theta\mu, \theta b)(\theta\upsilon)}
(\ass(\theta\mu, \theta b)) \bigr)\nonumber\\[-8pt]\\[-8pt]
&&\qquad =
\theta\bigl( \thetaa( \mu, b) ^{-1}
\bigl( \poiseed_{\thetaa( \mu, b)(\upsilon) }(\ass(\mu,
b)) \bigr) \bigr).\nonumber
\end{eqnarray}
To check (\ref{claimiso}), observe that by Lemma \ref{orient}
and the definition of the fixing isometry function,
\[
\thetaa(\theta\mu, \theta b) = \thetaa(\mu, b) \circ\theta
^{-1} .
\]
Hence, $\thetaa(\theta\mu, \theta b) ^{-1} = \theta
\circ\thetaa(\mu, b) ^{-1}$ and $\thetaa( \theta\mu, \theta
b)(\theta\upsilon) = \thetaa( \mu, b)(\upsilon)$. In addition, by Lemma
\ref{existenceassign}\hyperlink{assiso}{(b)}, $\ass(\mu, b) =
\ass(\theta\mu, \theta b)$, whence
\[
\poiseed_{\thetaa(
\theta\mu, \theta b)(\theta\upsilon)}(\ass(\theta\mu, \theta
b)) = \poiseed_{\thetaa( \mu, b)(\upsilon) }(\ass(\mu,
b)) .
\]
Thus, (\ref{claimiso}) holds.

Let $Y$ be a Poisson point process on $\R^d$ with
intensity $\ly> 0$. It follows from Lemmas
\ref{poiseed}, \ref{existenceassign} and the\vspace*{-1pt} fact that the
$d$-tags of all globes are nonsingular a.s. that
$\thmtwo(X) \eqd Y$, where $X$ is any Poisson point process on $\R^d$
with positive intensity. We need only note the following: the
Voronoi tessellation and the centers of the globes and Voronoi
cells depend only on $\Psi(X)$ (as in the case of $\phithmtwo$
from Proposition~\ref{phithmtwolabel}) and from Lemma
\ref{orient}, the fixing isometry function $\thetaa$ also
depends only on $(\Psi(X), X{\ronnn_{\Psi(X) ^c}})$.
\end{pf*}

Let us remark that the fact that the map $\thmtwo$ is
source-universal would not be very interesting without the additional
fact that it is strongly finitary, since using Theorem \ref{resulttwo}
we can define the following source-universal mapping. Let $\ly>0$.
For each $\lx>0$, let $\phi_{(\lx, \ly)}$ be the isometry-equivariant
mapping from Theorem \ref{resulttwo}, so that if $X$ is a Poisson point
process on $\R^d$ with intensity~$\lx$, then $\phi_{(\lx, \ly
)}(X)$ is
a Poisson point process on $\R^d$ with intensity $\ly$. Also, let
$J\dvtx\XX\to[0, \infty)$ be an isometry-invariant map such that for any
$\lx>0$, if $X$ is a Poisson point process on $\R^d$ with intensity
$\lx$, then $J(X) = \lx$ a.s. Clearly, the mapping
$\mathbf{\thmtwo} \dvtx\XX\to\XX$ defined by $\mu\mapsto\phi
_{(J(\mu),
\ly)}(\mu)$ is isometry-equivariant, and if $X$ is a Poisson point
process on $\R^d$ with positive intensity, then $\mathbf{\thmtwo}(X)$
is a Poisson point process on $\R^d$ with intensity~$\ly$.

\section{\texorpdfstring{Proof of Theorem \protect\ref{finitary}}{Proof of Theorem 4}}
\label{finitarysection}

In this section, we shall prove Theorem \ref{finitary} by
showing that the map $\Phii$ defined in (\ref{Phii}) and used
to prove Theorem \ref{result} and the map $\thmtwo$ defined
in (\ref{defthmtwo}) and used to prove Theorem \ref{resulttwo}
are both strongly finitary. We shall prove the following
stronger result from which Theorem \ref{finitary} follows
immediately.
\begin{theorem}
\label{finiteexp} Let $\Phii$ and $\thmtwo$ be the maps defined
in (\ref{Phii}) and (\ref{defthmtwo}), respectively. There
exists a map $T\dvtx\XX\to\N\cup\{ \infty \}$ such that if $X$
is a Poisson point process on $\R^d$ with positive intensity,
then $\E T(X)$ is finite and for all $\mu, \mu'\in\XX$ such
that $T(\mu) , T(\mu') < \infty$ and $\mu\ronnn_{B(\bfz,T)} =
\mu'\ronnn_{B(\bfz,T)}$, we have
\[
\Phii(\mu) \ronnn_{B(\bfz,1)} = \Phii(\mu') \ronnn_{B(\bfz,1)}
\quad\mbox{and}\quad \thmtwo(\mu) \ronnn_{B(\bfz,1)} = \thmtwo(\mu')
\ronnn_{B(\bfz,1)} .
\]
\end{theorem}

Let $\law$ be the law of $X$. Since $\E T(X) < \infty$ implies
that $T(X)$ is finite $\law$-a.s., Theorem \ref{finiteexp}
implies that $\Phii$ and $\thmtwo$ are both strongly finitary
with respect to $\law$.

We shall require the following additional property that the
selection rules defined in Section \ref{goodclumping}
satisfy.
\begin{lemma}
\label{local} Let $\Psi_{R}$ be the selection rule from Section
\ref{goodclumping}. For any $z \in\R^d$ any $\mu, \mu' \in
\XX$, if $\bar{B}(z,R)$ is a globe under $\mu$, then whenever
$\mu\ronnn_{B(z,5R+120+d)}=\mu'\ronnn_{B(z,5R+120+d)}$, we have
that $\bar{B}(z,R)$ is also a globe under $\mu'$.
\end{lemma}

Lemma \ref{local} is a localized version of property
\hyperlink{good}{(c)} in the definition of a~selection rule. We omit
the proof of Lemma \ref{local}, which uses the definition of
pre-seeds and seeds and is similar to that of Lemma
\ref{samepreseed}.
\begin{pf*}{Proof of Theorem \ref{finiteexp}}
Let $\Psi=\Psi_{R}$ be the $R$-selection rule from
Section~\ref{goodclumping} that is used to define the map $\Phii=\Phii_{(\lx,
\ly)}$. Recall that we use the $R$-selection rule with $R=1$ to define the
map $\thmtwo$. We now work toward a definition of $T$. Fix
$r:=100(5R + 101 +d)$. Let $\{ C_i \}_{i \in\Z^d}$ be an indexed partition
of $\R^d$ into\vspace*{1pt} equal-sized cubes of side length $r$ such that $C_i$ is
centered at $ir$. For all $i \in\Z^d$, let $c_i \subset C_i$ be the
ball of
radius $1$ concentric with the cube $C_i$, and let $E_i \subseteq\XX$
be the
set of measures such that $c_i$ contains a seed. Because the radius of $c_i$
is 1, it never contains more than one seed. Let $X$ be a Poisson point
process on $\R^d$ with positive intensity and law $\law$. It follows from
the definition of a seed and Lemma \ref{local} that $\{ {\mathbf1}_{X
\in
E_i} \}_{i \in\Z^d}$ is a collection of i.i.d. random variables
with positive
expectation. For each $i \in\Z^d$, let $E_i^1 \subset E_i$ be the set of
measures where the globe corresponding to the seed in $c_i$ is one-special,
and similarly let $E_i^2 \subset E_i$ be the set of measures where the globe
corresponding to the seed in $c_i$ is two-special. By
Proposition~\ref{fund}, it follows that $\{{\mathbf1}_{X \in E_i ^1}\}_{i \in\Z^d}$ and
$\{{\mathbf1}_{X \in E_i ^2}\}_{i\in\Z^d}$ are collections of i.i.d. random
variables with positive expectation. Let
\begin{eqnarray*}
T_{1}^{1}(\mu) &:=& \inf\bigl\{n \in\Z^{+} \dvtx\mu\in E^1_{(n, 0,
\ldots, 0)}
\mbox{ and for some } 0 < k_1 < k_2 < n, \\
&&\hspace*{95pt} \mbox{we
have } \mu
\in E^1_{(k_j, 0, \ldots, 0)} \mbox{ for } j=1,2 \bigr\}
\end{eqnarray*}
and
\begin{eqnarray*}
T_{-1}^{1}(\mu) &:=& \inf\bigl\{n \in\Z^{+} \dvtx\mu\in E^1_{(-n, 0,
\ldots, 0)}
\mbox{ and for some } 0 < k_1 < k_2 < n, \\
&&\hspace*{95.1pt} \mbox{we
have } \mu
\in E^1_{(-k_j, 0, \ldots, 0)} \mbox{ for } j=1,2 \bigr\} .
\end{eqnarray*}
Also define
\begin{eqnarray*}
T_1^{2}(\mu) &:=& \inf\bigl\{n \in\Z^{+} \dvtx\mu\in E^2_{(n, 0,
\ldots, 0)}
\mbox{ and for some } 0 < k_1 < k_2 < n, \\
&&\hspace*{95.4pt}
\mbox{we have }\mu\in
E^2_{(k_j, 0, \ldots, 0)}
\mbox{ for } j =1,2\bigr\}
\end{eqnarray*}
and
\begin{eqnarray*}
T_{-1}^{2}(\mu) &:=& \inf\bigl\{n \in\Z^{+} \dvtx\mu\in E^2_{(-n, 0,
\ldots, 0)}
\mbox{ and for some } 0 < k_1 < k_2 < n, \\
&&\hspace*{95.4pt}
\mbox{we have }\mu\in
E^2_{(-k_j, 0, \ldots, 0)}
\mbox{ for } j =1,2\bigr\} .
\end{eqnarray*}

Note that if we wanted to prove Theorem \ref{finiteexp} only for the
map $\Phii$, it would be enough to set $T= 8r(T_1^1 + T_1^2)$, but we
will require a slightly more complicated map $T$ to prove Theorem \ref
{finiteexp} for the map $\thmtwo$. Thus also similarly define
$T_i^{1}$, $T_i^{2}$, $T_{-i}^1$ and $T_{-i}^2$ for all $2 \leq i
\leq d$ by using coordinate $i$. Clearly, for all $ 1 \leq i
\leq d$, each of $T_i^{1}(X)$, $T_i^{2}(X)$, $T_{-i}^1(X)$ and
$T_{-i}^2(X)$ have finite mean.
We set
\[
T:= 8\sum_{i=1} ^d r(T_i^{1}+T_i^{2} + T_{-i}^1 + T_{-i}^2) .
\]

We now show that $\Phii$ satisfies the required property. Let
$\XX_T \subset\XX$ be the set of point measures such that $\mu
\in\XX_T$ iff $T(\mu) < \infty$. Observe that since $\Phii$
is monotone, to determine $\Phii(\mu) \ronnn_{B(\bfz,1)}$ it
suffices to determine which points of $[\mu] \cap
B(\bfz,1)$ will be in $[\Phii(\mu)] \cap B(\bfz,1)$. If
$x \in[\mu]$ does not belong to a globe, then whether or
not it is deleted depends on the value of $\ass(\mu, x)$.
Recall that $\ass$ is the assignment function for $\Psi$. If
$x \in[\mu] \cap B(\bfz,1)$ does belong to a globe, then
whether or not it is deleted depends on the globe $b$ for which
$x \in b$, on $\ass(\mu, b)$, and on the splitting
$\gammathin_b(\mu\ronnn_b, \ass(\mu, b))$. Let $\cglobe$ be
defined as in (\ref{centerglobe}), the point process of the
centers of the globes. Thus it suffices to show that for all
$\mu, \mu'\in\XX_T$ such that $\mu\ronnn_{B(\bfz,T(\mu))} =
\mu'\ronnn_{B(\bfz,T(\mu))}$, we have:
\begin{longlist}[(a)]
\item[(a)]
\hypertarget{sameglobes}
$\Psi(\mu) \cap B(\bfz,1) = \Psi(\mu') \cap B(\bfz,1)$ and
$\cglobe(\Psi(\mu)) \ronn_{B(\bfz,1)} =
\cglobe(\Psi(\mu')) \ronn_{B(\bfz,1)}$;
\item[(b)]
\hypertarget{sameass}
$\ass(\mu, x) = \ass(\mu', x)$ for all $x \in B(\bfz,1)$;
\item[(c)]
\hypertarget{sameassb}
$\ass(\mu, \bar{B}(y, R)) = \ass(\mu', \bar{B}(y, R))$ for all $y
\in B(\bfz,1)$.
\end{longlist}

Write $T = T(\mu)$.
Property \hyperlink{sameglobes}{(a)} follows from $\mu
\ronnn_{B(\bfz,T)} = \mu'\ronnn_{B(\bfz,T)}$ and
Lem\-ma~\ref{local}. Property \hyperlink{sameass}{(b)} follows from the
following observations. If $\bar{B}(z,R)$ is a~partner of some
$x \in B(\bfz,1)$, then $\|x-z\| \leq r\max(T_1^1,T_1^2)$.
Thus $B(z, 5R+120 +d) \subset B(\bfz,T)$. In addition, if $y
\in[\mu]$ shares a partner with $x$ and has lower one-
or two-rank than $x$, then $\|x-y\| \leq2r\max(T_1^1, T_1^2)$,
so that $y\in B(\bfz,T)$. Also note that the tag of a globe
$B(z,R)$ is contained in $B(z, 5R+120+d)$. Hence by
Lem\-ma~\ref{local}, the partners and ranks of $x$ are determined on
$B(\bfz,T)$. Thus
\[
p(\mu, x) = p(\mu', x)
\quad\mbox{and}\quad r(\mu, x) = r(\mu', x)
\]
for all $x \in B(\bfz,1)$ and all $\mu, \mu'\in\XX_T$ such
that $\mu\ronnn_{B(\bfz,T)} = \mu'\ronnn_{B(\bfz,T)}$, where
$p,r$ are the partner and rank functions of $\Psi$. From the
definition of $\ass$, property~\hyperlink{sameass}{(b)} follows.

Similarly, if $y \in B(\bfz,1)$, $b=\bar{B}(y,R)$ is a globe and
$\bar{B}(z,R)$ is a partner of~$b$, then $B(z, 5R+120+d)
\subset B(\bfz,T)$. In addition, if $\bar{B}(y',R)$ shares a~partner with $b$ and has a lower rank than $b$, then $B(y',5R+120+d)
\subset B(\bfz,T)$. Thus property
\hyperlink{sameassb}{(c)} also holds.

The proof that $\thmtwo$ has the required property is similar.
Recall that $\thmtwo$ is defined by placing Poisson point
processes inside each member of $\mathcal{V}(\cglobe)$ using
the assignment function $\ass$. Recall that $\mathcal{V}(\mu)$
is the Voronoi tessellation of the point process $\mu$, and each
Voronoi cell receives the $\uniforma$ variable assigned to the
globe that is contained in the Voronoi cell. For all $x \in
\R^d$, let $v(x, \mu)$ be the member of
$\mathcal{V}(\cglobe(\mu))$ to which $x$ belongs. From the
definition of $T$ and Lemma \ref{local}, it follows that for
all $x \in B(\bfz,1)$ and all $\mu, \mu'\in\XX_T$ such that
$\mu\ronnn_{B(\bfz,T(\mu))} = \mu'\ronnn_{B(\bfz,T(\mu))}$, we
have $v(x, \mu) = v(x, \mu')$. Moreover, it is not difficult
to verify that for each $x \in B(\bfz,1)$, if $b \subset v(x,
\mu)$ is a globe, then its partners, rank and assignment
function are also determined on $B(\bfz,T(\mu))$.
\end{pf*}

\section{Open problems}
\label{openprob}

Question \ref{thick}, in the \hyperref[intro]{Introduction}, asked whether
a~homogeneous Poisson point process $X$ on $\R^d$ can be
deterministically ``thickened'' via a factor---that is,
whether there exists a deterministic isometry-equivariant map
$\phi$ such that $\phi(X)$ is a homogeneous Poisson process of
higher intensity that contains all the original points of $X$.

We also do not know the answer to the following question, where
we drop the requirement of equivariance.
\begin{question}
\label{thicktwo} Let $d \geq1$, and let $\ly> \lx> 0$. Does
there exist a deterministic map $\phi$ such that if $X$ is a
homogeneous Poisson point process with intensity $\lx$, then $\phi
(X)$ is a
homogeneous Poisson point process on $\R^d$ with intensity
$\ly$, and such that \textit{all} the points of $X$ are points of
$\phi(X)?$
\end{question}

We can also ask similar questions in the discrete setting of
Bernoulli processes. We do not know the answer to following
simple question.
\begin{question}
\label{bernthick} Let $X =\{ X_i \}_{i \in\Z}$ be a sequence of
i.i.d. $\{ 0,1 \}$-valued random variables with $\E(X_0) =
\frac{1}{4}$. Does there exist a deterministic map $f$ such
that $\{ f(X)_i \}_{i \in\Z}$ is a sequence\vspace*{1pt} of i.i.d.
$\{ 0,1 \}$-valued random variables with $\E(f(X)_0) =
\frac{1}{2}$ and $f(X)_i \geq X_i$ for all $i \in\Z?$
\end{question}

Note that there does not exist a translation-equivariant map
$\phi$, a factor, that satisfies the condition of Question
\ref{bernthick}; if $\phi$ is a factor, then by the
Kolmogorov--Sinai theorem, the entropy of $\phi(X)$ cannot be
greater than the entropy of $X$. See \cite{MR1073173}, Chapter
5, for more details. More generally, if~$B(p)$
and $B(q)$ are Bernoulli shifts on $\{ 0,1,\ldots, d-1 \}$,
where the entropy of $p$ is less than the entropy of $q$, one
can ask whether there exists a deterministic map~$\phi$ from
$B(p)$ to $B(q)$ such that we have $\phi(x)_i \geq x_i$ for all
$x \in\{ 0, 1,\ldots, d-1 \}$ and all $i \in\Z$. Also see
\cite{balltwo} for more open problems.

\subsection*{Remark}
Ori Gurel-Gurevich and Ron Peled have informed us that they
have answered Questions \ref{thick}--\ref{bernthick} (with respective answers no, yes
and yes) in a~manuscript entitled ``Poisson Thickening'' \cite{GGP}.

\section*{Acknowledgments}
We thank the referee for providing a detailed and insightful report.
In particular, the referee pointed out a small error in an earlier
version of this paper and provided a somewhat simpler way to define a
seed given an equivalence class of pre-seeds.


%
%

%
\printaddresses


\begin{thebibliography}{31}

\bibitem{AHS}
%
\begin{barticle}[mr]
\bauthor{\bsnm{Angel},~\bfnm{Omer}\binits{O.}},
\bauthor{\bsnm{Holroyd},~\bfnm{Alexander~E.}\binits{A.~E.}} \AND
\bauthor{\bsnm{Soo},~\bfnm{Terry}\binits{T.}}
(\byear{2011}).
\btitle{Deterministic thinning of finite {P}oisson processes}.
\bjournal{Proc. Amer. Math. Soc.}
\bvolume{139}
\bpages{707--720}.
\bid{doi={10.1090/S0002-9939-2010-10535-X}, issn={0002-9939}, mr={2736350}}
\end{barticle}
%
\endbibitem

\bibitem{balltwo}
%
\begin{barticle}[mr]
\bauthor{\bsnm{Ball},~\bfnm{Karen}\binits{K.}}
(\byear{2005}).
\btitle{Monotone factors of i.i.d. processes}.
\bjournal{Israel J. Math.}
\bvolume{150}
\bpages{205--227}.
\bid{doi={10.1007/BF02762380}, issn={0021-2172}, mr={2255808}}
\end{barticle}
%
\endbibitem

\bibitem{ball}
%
\begin{barticle}[mr]
\bauthor{\bsnm{Ball},~\bfnm{Karen}\binits{K.}}
(\byear{2005}).
\btitle{Poisson thinning by monotone factors}.
\bjournal{Electron. Commun. Probab.}
\bvolume{10}
\bpages{60--69 (electronic)}.
\bid{issn={1083-589X}, mr={2133893}}
\end{barticle}
%
\endbibitem

\bibitem{Evans}
%
\begin{bincollection}[mr]
\bauthor{\bsnm{Evans},~\bfnm{Steven~N.}\binits{S.~N.}}
(\byear{2010}).
\btitle{A zero--one law for linear transformations of {L}\'evy noise}.
In \bbooktitle{Algebraic Methods in Statistics and Probability {II}}
(\beditor{M.~A. Viana} and \beditor{H.~P. Wynn}, eds.).
\bseries{Contemporary Mathematics}
\bvolume{516}
\bpages{189--197}.
\bpublisher{Amer. Math. Soc.}, \baddress{Providence, RI}.
\bid{mr={2730749}}
\bptnote{check year}%
\end{bincollection}
%
\endbibitem

\bibitem{MR2044812}
%
\begin{barticle}[mr]
\bauthor{\bsnm{Ferrari},~\bfnm{P.~A.}\binits{P.~A.}},
\bauthor{\bsnm{Landim},~\bfnm{C.}\binits{C.}} \AND
\bauthor{\bsnm{Thorisson},~\bfnm{H.}\binits{H.}}
(\byear{2004}).
\btitle{Poisson trees, succession lines and coalescing random walks}.
\bjournal{Ann. Inst. Henri Poincar\'e Probab. Stat.}
\bvolume{40}
\bpages{141--152}.
\bid{doi={10.1016/S0246-0203(03)00066-9}, issn={0246-0203}, mr={2044812}}
\end{barticle}
%
\endbibitem

\bibitem{GGP}
%
\begin{bmisc}[auto:STB|2011-03-03|12:04:44]
\bauthor{\bsnm{Gurel-Gurevich},~\bfnm{O.}\binits{O.}} \AND
\bauthor{\bsnm{Peled},~\bfnm{R.}\binits{R.}}
(\byear{2011}).
\bhowpublished{Poisson thickening. \textit{Israel J. Math.}
To appear. Available at
\href{http://arxiv.org/abs/arXiv:0911.5377}{arXiv:0911.5377}}.
\end{bmisc}
%
\endbibitem

\bibitem{random}
%
\begin{barticle}[mr]
\bauthor{\bsnm{Holroyd},~\bfnm{Alexander~E.}\binits{A.~E.}},
\bauthor{\bsnm{Pemantle},~\bfnm{Robin}\binits{R.}},
\bauthor{\bsnm{Peres},~\bfnm{Yuval}\binits{Y.}} \AND
\bauthor{\bsnm{Schramm},~\bfnm{Oded}\binits{O.}}
(\byear{2009}).
\btitle{Poisson matching}.
\bjournal{Ann. Inst. Henri Poincar\'e Probab. Stat.}
\bvolume{45}
\bpages{266--287}.
\bid{doi={10.1214/08-AIHP170}, issn={0246-0203}, mr={2500239}}
\end{barticle}
%
\endbibitem

\bibitem{MR1961286}
%
\begin{barticle}[mr]
\bauthor{\bsnm{Holroyd},~\bfnm{Alexander~E.}\binits{A.~E.}} \AND
\bauthor{\bsnm{Peres},~\bfnm{Yuval}\binits{Y.}}
(\byear{2003}).
\btitle{Trees and matchings from point processes}.
\bjournal{Electron. Commun. Probab.}
\bvolume{8}
\bpages{17--27 (electronic)}.
\bid{issn={1083-589X}, mr={1961286}}
\end{barticle}
%
\endbibitem

\bibitem{Extra-Heads}
%
\begin{barticle}[mr]
\bauthor{\bsnm{Holroyd},~\bfnm{Alexander~E.}\binits{A.~E.}} \AND
\bauthor{\bsnm{Peres},~\bfnm{Yuval}\binits{Y.}}
(\byear{2005}).
\btitle{Extra heads and invariant allocations}.
\bjournal{Ann. Probab.}
\bvolume{33}
\bpages{31--52}.
\bid{doi={10.1214/009117904000000603}, issn={0091-1798}, mr={2118858}}
\end{barticle}
%
\endbibitem

\bibitem{MR1084815}
%
\begin{bbook}[mr]
\bauthor{\bsnm{Horn},~\bfnm{Roger~A.}\binits{R.~A.}} \AND
\bauthor{\bsnm{Johnson},~\bfnm{Charles~R.}\binits{C.~R.}}
(\byear{1990}).
\btitle{Matrix Analysis}.
\bpublisher{Cambridge Univ. Press}, \baddress{Cambridge}.
\bnote{Corrected reprint of the 1985 original}.
\bid{mr={1084815}}
\end{bbook}
%
\endbibitem

\bibitem{jacod}
%
\begin{barticle}[mr]
\bauthor{\bsnm{Jacod},~\bfnm{Jean}\binits{J.}}
(\byear{1975}).
\btitle{Two dependent {P}oisson processes whose sum is still a {P}oisson
process}.
\bjournal{J.~Appl. Probab.}
\bvolume{12}
\bpages{170--172}.
\bid{issn={0021-9002}, mr={0370750}}
\end{barticle}
%
\endbibitem

\bibitem{MR1876169}
%
\begin{bbook}[mr]
\bauthor{\bsnm{Kallenberg},~\bfnm{Olav}\binits{O.}}
(\byear{2002}).
\btitle{Foundations of Modern Probability},
\bedition{2nd} ed.
\bpublisher{Springer}, \baddress{New York}.
\bid{mr={1876169}}
\end{bbook}
%
\endbibitem

\bibitem{keanea}
%
\begin{barticle}[mr]
\bauthor{\bsnm{Keane},~\bfnm{M.}\binits{M.}} \AND
\bauthor{\bsnm{Smorodinsky},~\bfnm{M.}\binits{M.}}
(\byear{1977}).
\btitle{A class of finitary codes}.
\bjournal{Israel J. Math.}
\bvolume{26}
\bpages{352--371}.
\bid{issn={0021-2172}, mr={0450514}}
\end{barticle}
%
\endbibitem

\bibitem{keaneb}
%
\begin{barticle}[mr]
\bauthor{\bsnm{Keane},~\bfnm{Michael}\binits{M.}} \AND
\bauthor{\bsnm{Smorodinsky},~\bfnm{Meir}\binits{M.}}
(\byear{1979}).
\btitle{Bernoulli schemes of the same entropy are finitarily isomorphic}.
\bjournal{Ann. of Math. (2)}
\bvolume{109}
\bpages{397--406}.
\bid{doi={10.2307/1971117}, issn={0003-486X}, mr={0528969}}
\end{barticle}
%
\endbibitem

\bibitem{kingman}
%
\begin{bbook}[mr]
\bauthor{\bsnm{Kingman},~\bfnm{J.~F.~C.}\binits{J.~F.~C.}}
(\byear{1993}).
\btitle{Poisson Processes}.
\bseries{Oxford Studies in Probability}
\bvolume{3}.
\bpublisher{Oxford Univ. Press}, \baddress{New York}.
\bid{mr={1207584}}
\end{bbook}
%
\endbibitem

\bibitem{last}
%
\begin{barticle}[mr]
\bauthor{\bsnm{Last},~\bfnm{G{\"u}nter}\binits{G.}} \AND
\bauthor{\bsnm{Thorisson},~\bfnm{Hermann}\binits{H.}}
(\byear{2009}).
\btitle{Invariant transports of stationary random measures and
mass-stationarity}.
\bjournal{Ann. Probab.}
\bvolume{37}
\bpages{790--813}.
\bid{doi={10.1214/08-AOP420}, issn={0091-1798}, mr={2510024}}
\end{barticle}
%
\endbibitem

\bibitem{MR2132405}
%
\begin{bbook}[mr]
\bauthor{\bsnm{Molchanov},~\bfnm{Ilya}\binits{I.}}
(\byear{2005}).
\btitle{Theory of Random Sets}.
\bpublisher{Springer}, \baddress{London}.
\bid{mr={2132405}}
\end{bbook}
%
\endbibitem

\bibitem{MR0447525}
%
\begin{bbook}[mr]
\bauthor{\bsnm{Ornstein},~\bfnm{Donald~S.}\binits{D.~S.}}
(\byear{1974}).
\btitle{Ergodic Theory, Randomness, and Dynamical Systems}.
\bpublisher{Yale Univ. Press}, \baddress{New Haven}.
\bid{mr={0447525}}
\end{bbook}
%
\endbibitem

\bibitem{MR910005}
%
\begin{barticle}[mr]
\bauthor{\bsnm{Ornstein},~\bfnm{Donald~S.}\binits{D.~S.}} \AND
\bauthor{\bsnm{Weiss},~\bfnm{Benjamin}\binits{B.}}
(\byear{1987}).
\btitle{Entropy and isomorphism theorems for actions of amenable groups}.
\bjournal{J. Anal. Math.}
\bvolume{48}
\bpages{1--141}.
\bid{issn={0021-7670}, mr={0910005}}
\end{barticle}
%
\endbibitem

\bibitem{MR1073173}
%
\begin{bbook}[mr]
\bauthor{\bsnm{Petersen},~\bfnm{Karl}\binits{K.}}
(\byear{1989}).
\btitle{Ergodic Theory}.
\bseries{Cambridge Studies in Advanced Mathematics}
\bvolume{2}.
\bpublisher{Cambridge Univ. Press}, \baddress{Cambridge}.
\bnote{Corrected reprint of the 1983 original}.
\bid{mr={1073173}}
\end{bbook}
%
\endbibitem

\bibitem{MR681482}
%
\begin{bbook}[mr]
\bauthor{\bsnm{Rees},~\bfnm{Elmer~G.}\binits{E.~G.}}
(\byear{1983}).
\btitle{Notes on Geometry}.
\bpublisher{Springer}, \baddress{Berlin}.
\bid{mr={0681482}}
\end{bbook}
%
\endbibitem

\bibitem{MR1199815}
%
\begin{bbook}[mr]
\bauthor{\bsnm{Reiss},~\bfnm{R.~D.}\binits{R.~D.}}
(\byear{1993}).
\btitle{A Course on Point Processes}.
\bpublisher{Springer}, \baddress{New York}.
\bid{mr={1199815}}
\end{bbook}
%
\endbibitem

\bibitem{MR2306207}
%
\begin{bincollection}[mr]
\bauthor{\bsnm{Serafin},~\bfnm{Jacek}\binits{J.}}
(\byear{2006}).
\btitle{Finitary codes, a short survey}.
In \bbooktitle{Dynamics \& Stochastics}.
\bseries{Institute of Mathematical Statistics Lecture Notes---Monograph Series}
\bvolume{48}
\bpages{262--273}.
\bpublisher{IMS}, \baddress{Beachwood, OH}.
\bid{mr={2306207}}
\end{bincollection}
%
\endbibitem

\bibitem{MR0161960}
%
\begin{barticle}[mr]
\bauthor{\bsnm{Sina{\u\i}},~\bfnm{Ja.~G.}\binits{J.~G.}}
(\byear{1962}).
\btitle{A weak isomorphism of transformations with invariant measure}.
\bjournal{Dokl. Akad. Nauk SSSR}
\bvolume{147}
\bpages{797--800}.
\bid{issn={0002-3264}, mr={0161960}}
\end{barticle}
%
\endbibitem

\bibitem{soo-2006}
%
\begin{barticle}[mr]
\bauthor{\bsnm{Soo},~\bfnm{Terry}\binits{T.}}
(\byear{2010}).
\btitle{Translation-equivariant matchings of coin flips on {$\Bbb Z\sp d$}}.
\bjournal{Adv. in Appl. Probab.}
\bvolume{42}
\bpages{69--82}.
\bid{doi={10.1239/aap/1269611144}, issn={0001-8678}, mr={2666919}}
\end{barticle}
%
\endbibitem

\bibitem{borel}
%
\begin{bbook}[mr]
\bauthor{\bsnm{Srivastava},~\bfnm{S.~M.}\binits{S.~M.}}
(\byear{1998}).
\btitle{A Course on {B}orel Sets}.
\bseries{Graduate Texts in Mathematics}
\bvolume{180}.
\bpublisher{Springer}, \baddress{New York}.
\bid{mr={1619545}}
\end{bbook}
%
\endbibitem

\bibitem{Thorissontrv}
%
\begin{barticle}[mr]
\bauthor{\bsnm{Thorisson},~\bfnm{Hermann}\binits{H.}}
(\byear{1996}).
\btitle{Transforming random elements and shifting random fields}.
\bjournal{Ann. Probab.}
\bvolume{24}
\bpages{2057--2064}.
\bid{doi={10.1214/aop/1041903217}, issn={0091-1798}, mr={1415240}}
\end{barticle}
%
\endbibitem

\bibitem{MR1741181}
%
\begin{bbook}[mr]
\bauthor{\bsnm{Thorisson},~\bfnm{Hermann}\binits{H.}}
(\byear{2000}).
\btitle{Coupling, Stationarity, and Regeneration}.
\bpublisher{Springer}, \baddress{New York}.
\bid{mr={1741181}}
\end{bbook}
%
\endbibitem

\bibitem{timarb}
%
\begin{bmisc}[auto:STB|2011-03-03|12:04:44]
\bauthor{\bsnm{Tim{\'a}r},~\bfnm{{\'A}.}\binits{{\'A}.}}
\bhowpublished{Invariant matchings of exponential tail on coin flips in
$\mathbb{Z}^d$. Preprint. Available at \href{http://arxiv.org/abs/arXiv:0909.1090}{arXiv:0909.1090}}.
\end{bmisc}
%
\endbibitem

\bibitem{MR2081459}
%
\begin{barticle}[mr]
\bauthor{\bsnm{Tim{\'a}r},~\bfnm{{\'A}d{\'a}m}\binits{{\'A}.}}
(\byear{2004}).
\btitle{Tree and grid factors for general point processes}.
\bjournal{Electron. Commun. Probab.}
\bvolume{9}
\bpages{53--59 (electronic)}.
\bid{issn={1083-589X}, mr={2081459}}
\end{barticle}
%
\endbibitem

\end{thebibliography}
\end{document}